# LOWER BOUNDS FOR THE DENSITY OF LOCALLY ELLIPTIC ITÔ PROCESSES

By Vlad Bally

*Université de Marne-la-Vallée*

We give lower bounds for the density $p_T(x,y)$ of the law of $X_t$, the solution of $dX_t = \sigma(X_t)\,dB_t + b(X_t)\,dt, X_0 = x$, under the following local ellipticity hypothesis: there exists a deterministic differentiable curve $x_t$, $0 \leq t \leq T$, such that $x_0 = x$, $x_T = y$ and $\sigma\sigma^*(x_t) > 0$, for all $t \in [0,T]$. The lower bound is expressed in terms of a distance related to the skeleton of the diffusion process. This distance appears when we optimize over all the curves which verify the above ellipticity assumption.

The arguments which lead to the above result work in a general context which includes a large class of Wiener functionals, for example, Itô processes. Our starting point is work of Kohatsu-Higa which presents a general framework including stochastic PDE's.

**1. Introduction.** It is well known that under uniform ellipticity and boundedness assumptions for the diffusion coefficients matrix, the law of a diffusion process is absolutely continuous with respect to the Lebesgue measure and one may obtain Gaussian-type lower and upper bounds for the density of the law. This classical result has been extended (see [7, 12, 17]) to the more subtle case where, instead of ellipticity, one assumes a Hörmander-type hypothesis. In this paper, we do not proceed in this direction. On the other hand, as an application of Malliavin's calculus, it is proven that under appropriate hypothesis, a large variety of functionals on the Wiener space (e.g., solutions of stochastic PDE's) have absolute continuous laws and the density is smooth (see [16]). Using already standard techniques, one may prove that some Gaussian upper bounds hold true. In a number of cases, one may also succeed to prove that the density is strictly positive (see e.g., [1, 3, 5, 15] or [16]). But the techniques used to prove strict positivity are rather qualitative and do not provide lower bounds. So, this remains a challenging problem. In a recent paper, Kohatsu-Higa [13] developed a strategy which permits an attack on this problem for abstract Wiener functionals. The author proposes a framework which essentially expresses the idea of uniform ellipticity for a Wiener functional and then develops a methodology for computing lower bounds. He employs this method for the stochastic heat equation. More recently, Dalang and Nualart [6] provided applications to potential theory for hyperbolic SPDE's.

The paper of Kohatsu-Higa was the starting point for our work and several important ideas come from it. But we give a local approach which permits the treatment of a significantly larger class of problems. On one hand, we avoid boundedness assumptions on the coefficients of the equations at hand. In recent work, Guérin, Méléard and Nualart [9] used this local approach in order to obtain lower bounds for the solution of Landau's equation—a serious difficulty there is that the coefficients are not bounded. But the main purpose is to relax the uniform ellipticity hypothesis: we simply assume that there exists a deterministic differentiable curve such that the ellipticity assumption holds true along this curve. This gives access to a large class of problems which are far from uniform elliptic diffusions, such as stochastic integrals and solutions of non-Markov stochastic equations (see the examples in [2]). These problems are also out of reach of the criterion based on Hörmander's hypothesis (but the method presented here does not cover this criterion).

Although our main applications concern diffusion processes, we present the method in a more general context which is close to the abstract setting put forward by Kohatsu-Higa. We consider a $q$-dimensional Itô process of the form

$$X_t^i = x_0^i + \sum_{j=1}^{\infty} \int_0^t U_s^{ij} \, dB_s^j + \int_0^t V_s^i \, ds, \qquad i = 1, \ldots, q,$$

where $B^j, j \in N$, are independent Brownian motions. We are interested in the density $p_T(x_0, y)$ of $X_T$ at a point $y$. We assume that $U$ and $V$ are smooth in Malliavin's sense so that $X_T$ is also smooth. We now give the nondegeneracy assumption. We fix a deterministic differentiable curve $x_t, 0 \leq t \leq T$, such that $x_0 = x_0$, $x_T = y$ and some deterministic functions $r_t, K_t > 0$ for $0 \leq t \leq T$. We also consider a family of deterministic $q \times q$ symmetric positive definite matrices $Q_t, 0 \leq t \leq T$, and denote by $\lambda_t > 0$ the lower eigenvalue of $Q_t$. Given $t$ and $\delta > 0$, we define

$$\Gamma_\delta^i(t) := \sum_{j=1}^{\infty} \int_t^{t+\delta} (U_s - U_t)^{ij} \, dB_s^j + \int_t^{t+\delta} V_s^i \, ds.$$

Then our hypotheses are the following. For every $0 < t < T$ and $0 < \delta < T - t$,

(H, i)  $\qquad\qquad\qquad\qquad U_t U_t^* \geq Q_t,$

$(\mathrm{H}^\nu, \mathrm{ii})$ $\qquad \|\Gamma_\delta(t)\|_{k,p,t} \leq K(t)\delta^{1/2+\nu}, \qquad \nu > 0,$

on the set defined by $|Q_t^{-1/2}(X(t) - x(t))| \leq r(t)$.

Let us explain this definition. One writes

(1) $$X^i_{t+\delta} = X^i_t + \sum_{j=1}^{\infty} U^{ij}_t (B^j_{t+\delta} - B^j_t) + \Gamma^i_\delta(t).$$

The random variable $G_\delta(t) =: X_t + \sum_{j=1}^{\infty} U^j_t(B^j_{t+\delta} - B^j_t)$ is Gaussian conditionally to $F_t = \sigma(B^j_s, s \leq t, j \in N)$ and has the covariance matrix $\delta \times U_t U_t^*$. (H, i) therefore says that this term is nondegenerate; it represents the ellipticity assumption. $\Gamma^i_\delta(t)$ is a remainder and $(\mathrm{H}^\nu, \mathrm{ii})$ says that this remainder may be ignored with respect to the principal term $G_\delta(t)$, which is essentially of order $\delta^{1/2}$. $\nu$ is a strictly positive number which depends on the problem at hand—in the context of diffusion processes, $\nu = \frac{1}{2}$ and for the stochastic heat equation, $\nu = \frac{1}{4}$ (see [13]). The norm $\| \circ \|_{k,p,t}$ is a Sobolev norm which involves the $L^p$-norms of the first $k$ Malliavin derivatives where $p, k$ are some integers depending on the dimension $q$. The lower index $t$ signifies that we work with conditional expectations with respect to $F_t$ and not with usual expectations; we use a conditional version of the Malliavin calculus. Let us now comment on the localization. Both $U_t = U_t(\omega)$ and $\|\Gamma_\delta(t)\|_{k,p,t} = \|\Gamma_\delta(t)\|_{k,p,t}(\omega)$ are random variables. So, the hypotheses (H, i) and $(\mathrm{H}^\nu, \mathrm{ii})$ hold true only for $\omega \in \{|Q_t^{-1/2}(X(t) - x(t))| \leq r(t)\}$. Let us consider the example of the diffusion process $dX_t = \sigma(X_t) dB_t + b(X_t) dt$. Then $U_t = \sigma(X_t)$ and so (H, i) says that $\sigma\sigma^*(x) \geq Q_t$ for $x$ such that $|Q_t^{-1/2}(x - x(t))| \leq r(t)$. Therefore we need the ellipticity assumption only on a tube around the curve $x_t$.

Roughly speaking, in order to obtain lower bounds for $p_T(x_0, y)$, we proceed as follows. We construct a time grid $0 = t_0 < \cdots < t_N = T$ and let $\delta_i = t_i - t_{i-1}$. We denote by $p_i(z)$ [resp. $\overline{p}_i(z)$] the $F_{t_i}$-conditional density of $X(t_{i+1})$ [resp. of $G_{\delta_{i+1}}(t_i)$] with respect to the Lebesgue measure. We first note that if $|Q_{t_i}(X(t_i) - z)| \leq \delta_i$, then $\overline{p}_i(z) \geq 1/e^2(2\pi\delta_i)^{d/2}$. This is an easy computation based on the fact that $G_{\delta_{i+1}}(t_i)$ is a Gaussian random variable and we control the covariance matrix by means of (H, i). Next, we want to use the fact that the reminder $\Gamma_{\delta_{i+1}}(t_i)$ is small in order to derive a similar evaluation for $p_i(z)$. This is a more involved computation because $p_i(z) = E\delta_z(X(t_{i+1})) = E\delta_z(G_{\delta_{i+1}}(t_i) + \Gamma_{\delta_{i+1}}(t_i))$ where $\delta_z$ is the Dirac function. Since the Dirac function is not smooth, the fact that $\Gamma_{\delta_{i+1}}(t_i)$ is small in $L^p$-norms is not sufficient—we need the Sobolev norms (in Malliavin's sense) to also be small—this is why $\|\Gamma_{\delta_{i+1}}(t_i)\|_{k,p,t}$ appears in $(\mathrm{H}^\nu, \mathrm{ii})$. We may then use a development in Taylor series and Malliavin's integration

by parts formula (this is very similar to the calculus in [13], except for a localization argument which allows the avoidance of uniform ellipticity assumptions). This evaluation represents the basic element in the calculus and now our problem now is to transport it, by means of a "chain argument," along the curve $x_t$. This is done in the abstract context of the "evolution sequences" in Section 2. In Section 3, we discuss the Itô processes presented before and in Section 4, we deal with diffusion processes.

There is a certain analogy between the strategy used here and the one employed in the analytical approach to this problem (compare the decomposition used in (1) with (4.1), (4.2) page 14 in [8]; see also [4]). The advantage of the stochastic method is that it permits localization on the set of trajectories which remain in a tube around the deterministic curve. This allows the treatment of certain classes of diffusions which are not uniform elliptic and which do not have bounded coefficients. But the drawback is that we need much more regularity for the coefficients of the diffusion process.

In the context of diffusion processes, we are able to give a nice form of the lower bound by means of a distance based on the skeleton of the diffusion process. More precisely, suppose that $X_t \in R^q$, $t \geq 0$, solves the SDE

$$(*) \qquad dX_t = \sum_{j=1}^{d} \sigma_j(X_t)\, dB_t^j + b(X_t)\, dt, \qquad X_0 = x_0.$$

We denote $p_T(x_0, y) = P(X_T \in dy)$. We assume that the coefficients have linear growth, are $q+2$ times differentiable and have bounded derivatives. Moreover, we consider some functions $\lambda_*, \lambda^* : R^q \to R_+$ and assume that $\lambda^*(x) \geq \sigma\sigma^*(x) \geq \lambda_*(x) \geq 0$. In particular, $\lambda_*$ may be the lower eigenvalue of $\sigma\sigma^*$, but for technical reasons, we accept smaller functions as well. Finally, we consider a control $\phi = (\phi^1, \ldots, \phi^d)$, $\phi^j \in L^2[0, T]$, and denote by $x^\phi$ the solution of the ordinary differential equation

$$(**) \qquad dx_t^\phi = \sum_{j=1}^{d} \sigma_j(x_t^\phi) \phi_t^j dt, \qquad x_0^\phi = x_0.$$

We consider a set of parameters $\theta = (\mu, \chi, \nu, \eta, h)$, $\mu, \nu, \eta \geq 1$, $h, \chi > 0$, and we define $\Phi_\theta(x_0, y)$ to be the set of the controls $\phi \in (L^2([0,T]))^d$ such that

$$x_0^\phi = x_0, \qquad x_T^\phi = y, \qquad \frac{\lambda_*(x_t^\phi)}{\lambda^*(x_t^\phi)} \geq \frac{1}{\mu}, \qquad \sqrt{\lambda_*(x_t^\phi)} \geq \frac{1}{\chi} \qquad \forall t \in [0, T],$$

$$|\phi_t| \leq \eta |\phi_s| \qquad \forall\, |s - t| \leq h, \qquad |\phi_t| \leq \nu \qquad \forall t \leq T.$$

Then we define

$$d_\theta(x_0, y) = \inf\left\{ \|\phi\|_T = \left(\int_0^T |\phi_t|^2 dt\right)^{1/2} : \phi \in \Phi_\theta(x_0, y) \right\}$$

$$= \infty \qquad \text{if } \Phi_\theta(x_0, y) = \varnothing.$$

Our lower bound is given by

$$p_T(x_0, y) \geq \frac{1}{4e^2(6\mu\sqrt{q}\pi T)^{q/2}\sqrt{\det \sigma\sigma^*(y)}}$$
$$\times \exp\Big(-K_q(1 + \ln(\mu\eta))$$
$$\times \Big(\mu^4 d_\theta^2(x_0, y) + T\Big(\mu^4 \vee (\mu + \chi)^2 K_{\text{diff}} + \frac{1}{h} + \nu\Big)\Big)\Big),$$

where $K_{\text{diff}}$ depends on the bounds of the diffusion coefficients and $K_q$ is a constant depending only on $q$.

## 2. Evolution sequences.

2.1. *Conditional Malliavin calculus.* We consider a probability space $(\Omega, F, P)$ with a filtration $F_t$, $t \geq 0$, and an infinite-dimensional Brownian motion $B = (B^j)_{j \in N}$ with respect to this filtration (we do not need to assume that the filtration $F_t$ is generated by the Brownian motion itself). Moreover, we fix some $t \geq 0$, $\delta > 0$ and denote by $E_t$ the conditional expectation with respect to $F_t$, that is,

$$E_t(\Phi) := E(\Phi|F_t).$$

We will use a conditional version of Malliavin's calculus that we shall now outline. We work with the standard Malliavin derivative operators, but we will consider some specific norms which permit vs to focus on the derivatives with respect to $B_s, s \in [t, t+\delta]$ (instead of $s \in [0, \infty)$), on one hand and we will replace the expectation $E$ by the conditional expectation $E_t$ on the other hand. Let us briefly recall some notation (we refer to [16] or to [14] for a complete exposition of this topic). $D^{k,p}$ is the space of random variables which are $k$ times differentiable in Malliavin's sense, in $L^p$. For $F \in D^{k,p}$, the derivative of order $k$ is $D^k F$, an element of the space $H_k$ which is defined in the following way. We denote by $\Theta_k$ the set of the multi-indices $\alpha = (\alpha_1, \ldots, \alpha_k), \alpha_i \in \{1, 2, \ldots\}$ and let $R^{\Theta_k} = \{(x_\alpha)_{\alpha \in \Theta_k} : x_a \in R\}$. For a measurable function $V : [0, \infty)^k \to R^{\Theta_k}$, we define

$$|V|_k^2 := \int_{[0,\infty)^k} \sum_{\alpha \in \Theta_k} |V^\alpha(s_1, \ldots, s_k)|^2 \, ds_1, \ldots, ds_k \quad \text{and}$$

$$H_k := \{V : [0, \infty)^k \to R^{\Theta_k} : |V|_k^2 < \infty\}.$$

$H_k$ is a Hilbert space with the scalar product

$$\langle V, U \rangle_k := \int_{[0,\infty)^k} \sum_{\alpha \in \Theta_k} V^\alpha(s_1, \ldots, s_k) U^\alpha(s_1, \ldots, s_k) \, ds_1, \ldots, ds_k.$$

For $F \in D^{k,p}$, we denote by $D^k F$ the derivative of order $k$, that is, $D^k_{s_1,\ldots,s_k} F = (D^{k,\alpha}_{s_1,\ldots,s_k} F)_{\alpha \in \Theta_k}$ (see [16]) and we have $E|D^k F|^p_k < \infty$. So, $D^k F \in H_k$.

The above scalar product is used in the standard Malliavin calculus. In our framework, for every fixed $t, \delta > 0$ we define

$$\langle V, U \rangle_{t,\delta,k} := \int_{[t,t+\delta)^k} \sum_{\alpha \in \Theta_k} V^\alpha(s_1,\ldots,s_k) U^\alpha(s_1,\ldots,s_k)\, ds_1,\ldots,ds_k,$$

$$|V|^2_{t,\delta,k} := \langle V, V \rangle_{t,\delta,k} = \int_{[t,t+\delta)^k} \sum_{\alpha \in \Theta_k} |V^\alpha(s_1,\ldots,s_k)|^2\, ds_1,\ldots,ds_k.$$

For $F \in D^{k,p}$, we define the following Sobolev norms:

$$\|F\|^2_{t,\delta,k} := \sum_{i=0}^k |D^i F|^2_{t,\delta,i} = \sum_{i=0}^k \sum_{\alpha \in \Theta_i} \int_{[t,t+\delta)^i} |D^{i,\alpha}_{s_1,\ldots,s_i} F|^2\, ds_1,\ldots,ds_i,$$

$$\|F\|^p_{t,\delta,k,p} := E_t(\|F\|^p_{t,\delta,k}), \qquad \|\!|F|\!\|^p_{t,\delta,k,p} := \|F\|^p_{t,\delta,k,p} - E_t(|F|^p).$$

Moreover, for a multidimensional functional $F = (F_1,\ldots,F_q)$, we denote $\|F\|^p_{t,\delta,k,p} = \sum_{i=1}^q \|F_i\|^p_{t,\delta,k,p}$.

Notice that $\|F\|_{t,\delta,k,p}$ is not a constant (as in the standard case), but an $F_t$-measurable random variable. Notice also that using $\|F\|_{t,\delta,k,p}$ instead of the standard norm $\|F\|_{k,p} = \|F\|_{0,\infty,k,p}$ amounts to using the differential calculus with respect to $B_s, s \in [t, t+\delta]$, only and taking conditional expectations instead of the usual expectation.

In the sequel, we will employ the following inequality. There exists a universal constant $\mu(k)$ such that for every $F \in D^{k,p}$, every $p > 1$ and every $\phi \in C_b^k$, one has

$$\|\!|\phi(F)|\!\|_{t,\delta,k,p} \leq \mu(k) \|\phi\|_{k,\infty} \|\!|F|\!\|_{t,\delta,k,2^k p}, \tag{2}$$

where $\|\phi\|_{k,\infty} := \max_{i=0,k} \sup_{x \in R} |\phi^{(i)}(x)|$. This is a straightforward consequence of the chain rule and inequality (22) from the Appendix.

We now define the Ornstein–Uhlembeck operator $L_{t,\delta}$ associated to $B_s, s \in [t, t+\delta]$, by

$$L_{t,\delta} F = \sum_{i=1}^\infty \int_t^{t+\delta} D^{1,i}_s F\, dB^i_s,$$

where the above stochastic integral is the Skorohod integral (see [16] or [14]). If $F \in \mathrm{Dom}(L)$, where $L$ is the standard Ornstein–Uhlembeck operator, then $L_{t,\delta} F$ is well defined. In our framework, we will assume that $F \in \bigcap_{p \in N} D^{k+2,p}$ for some $k \in N$ and so, in view of Meyer's inequalities, $F \in \mathrm{Dom}(L)$ and $LF \in \bigcap_{p \in N} D^{k,p}$. We will use the following form of Meyer's inequality which is proved in [16]: there exists an universal constant $c(k,p)$ such that for every $F \in \bigcap_{r \in N} D^{k+2,r}$,

$$\|L_{t,\delta}F\|_{t,\delta,k,p} \leq c(k,p)\|F\|_{t,\delta,k+2,p}. \tag{3}$$

In the sequel, we assume that $c(k,p)$ increases in both $k$ and $l$. If not, we take the maximum over $k' \leq k$ and $p' \leq p$.

It is easy to check (the standard argument) that for $F, G \in D^{2,2}$,

$$E_t(\langle DF, DG \rangle_{t,\delta,1}) = E_t(FL_{t,\delta}G) = E_t(GL_{t,\delta}F).$$

Here and in the sequel, we use the notation $D$ instead of $D^1$.

This is a conditional version of the standard duality relation which is the starting point for Malliavin's calculus. The same arguments as those used in the classical case give the conditional version of the integration by parts formula presented in the following theorem. Before stating this result, we define the Malliavin covariance matrix corresponding to $[t, t+\delta]$. Let $F = (F_1, \ldots, F_q), F_1, \ldots, F_q \in D^{1,2}$. We define

$$\phi_{t,\delta,F}^{ij} := \langle DF_i, DF_j \rangle_{t,\delta,1}, \qquad i,j = 1, \ldots, q.$$

We now state a localized version of the nondegeneracy assumption in Malliavin calculus. We consider a measurable set $A \subset \{\det \phi_{t,\delta,F} \neq 0\}$ and denote

$$v_p(F, A) = (E_t(1 \vee (\det \phi_{t,\delta,F})^{-p} \mathbb{1}_A))^{1/p}.$$

We assume that $v_p(F, A) < \infty$ for every $p \in N$ and denote by $\widehat{\phi}_{t,\delta,F}(\omega)$ the inverse of $\phi_{t,\delta,F}(\omega)$ for $\omega \in A$. We denote by $D_A^k$ the class of random variables $G \in \bigcap_{p \in N} D^{k,p}$ such that $G(\omega) = 0$ and $D^i G(\omega) = 0$, $i = 1, \ldots, k$, for $\omega \in A^c$. The following lemma gives the localized version of the nondegenerancy condition:

LEMMA 1. *Let $F = (F_1, \ldots, F_q)$ with $F_i \in D_A^{k+1}$, $i = 1, \ldots, q$, and $G = (G_1, \ldots, G_q)$ with $G_i \in D_A^k$, $i = 1, \ldots, q$. Assume that $v_p(F, A) < \infty$ for every $p \in N$. Then $\widehat{\phi}_{t,\delta,F} \times G \in (\bigcap_{p \in N} D^{k,p})^q$. Moreover, there exists a universal constant $c(k,q)$ such that for every $p \geq 1$,*

$$\|(\widehat{\phi}_{t,\delta,F}G)^i\|_{t,\delta,k,p} \leq c(k,q)\|G\|_{t,\delta,k,4p}v_{8(k+1)p}^{k+1}(F,A)\|F\|_{t,\delta,k+1,2^{2(q+2)}p}^{4q+2}. \tag{4}$$

The proof is straightforward and so we leave it for the Appendix.

The same proof as in the standard case gives the following integration by parts theorem:

THEOREM 2. *Let $F = (F_1, \ldots, F_q)$ with $F_1, \ldots, F_q \in \bigcap_{p \in N} D^{2,p}$. Consider a measurable set $A$ such that $v_p(F, A) < \infty$ for every $p \in N$ and a* *random variable* $G \in D_A^1$. *Then for every smooth function* $f \colon R^q \to R$ *and every* $i = 1, \ldots, q$, *one has*

$$(\text{IP}_i) \quad E_t\left(\frac{\partial f}{\partial x_i}(F)G\right) = E_t(f(F)H_i(F,G))$$

$$\text{with } H_i(F,G) = -\sum_{j=1}^{q}(G\widehat{\phi}_{t,\delta,F}^{ji}L_{t,\delta}(F_j) + \langle DF_j, D(\widehat{\phi}_{t,\delta,F}^{ji}G)\rangle_{t,\delta,1}).$$

*Suppose that* $F_1, \ldots, F_q \in \bigcap_{p \in N} D^{k+1,p}$ *and* $G \in D_A^k$ *for some* $k \in N$. *Then for every multi-index* $\alpha = (\alpha_1, \ldots, \alpha_k) \in \{1, \ldots, q\}^k$, *one has*

$$(\text{IP}_\alpha) \quad \begin{aligned} E_t(D^\alpha f(F)G) &= E_t(f(F)H_\alpha(F,G)) \\ \text{with } H_\alpha(F,G) &= H_{\alpha_k}(F, H_{(\alpha_1, \ldots, \alpha_{k-1})}(F,G)), \end{aligned}$$

*where* $D^\alpha f = \partial^k f / \partial x_{\alpha_1}, \ldots, \partial x_{\alpha_k}$.

We now give some evaluations of the norms of $H_\alpha(F,G)$.

PROPOSITION 3. *Suppose that* $F_1, \ldots, F_q \in \bigcap_{p \in N} D^{k+l+2,p}$ *and* $G \in D_A^{l+1}$ *for some* $k, l \in N$. *Then there exists a universal constant* $c(k, l, q)$ *such that for every multi-index* $\alpha$ *with* $|\alpha| = l$ *and every* $p \in N$, *one has*

$$(5) \quad \|H_\alpha(F,G)\|_{t,\delta,k,p} \leq c(k,l,q) \times \|G\|_{t,\delta,k+l,2^{3l}p} v_{2^{4l}(k+l+1)^l p}^{l(k+l)}(F,A)\|\|F\|\|_{t,\delta,k+l+1,2^{l(q+4)}p}^{l(4q-1)}.$$

*In particular (taking* $k = 0$, $l = q+1$ *and* $p = 1$*), there exists two universal constants depending only on* $q$, $c^* = c^*(q)$ *and* $p^* = p^*(q)$, *such that*

$$(6) \quad E_t(|H_\alpha(F,G)|) \leq c^*\|G\|_{t,\delta,q+1,p^*} v_{p^*}^{(q+1)^2}(F,A)\|\|F\|\|_{t,\delta,q+2,p^*}^{(q+1)(4q-1)}$$

*for every multi-index* $\alpha$ *with* $|\alpha| \leq q+1$.

We leave the proof for the Appendix.

To finish, we give the following simple fact concerning the Malliavin covariance matrix. We denote by $\underline{\lambda}_{t,\delta,F}$ (resp. $\overline{\lambda}_{t,\delta,F}$) the smaller (resp. the larger) eigenvalue of $\phi_{t,\delta,F}$. They are given by

$$\underline{\lambda}_{t,\delta,F} = \inf_{|\xi|=1} \sum_{i,j=1}^{q} \xi_i \xi_j \phi_{t,\delta,F}^{ij}, \qquad \overline{\lambda}_{t,\delta,F} = \sup_{|\xi|=1} \sum_{i,j=1}^{q} \xi_i \xi_j \phi_{t,\delta,F}^{ij}.$$

PROPOSITION 4. *Let* $F, G \in (D^{1,2})^q$. *Then*

$$(\det \phi_{t,\delta,F+G})^{1/q} \geq \tfrac{1}{2}\underline{\lambda}_{t,\delta,F} - \overline{\lambda}_{t,\delta,G}.$$

PROOF. Using the elementary inequality $(x+y)^2 \geq \frac{1}{2}x^2 - y^2$, one obtains

$$(\det \phi_{t,\delta,F+G})^{1/q} \geq \underline{\lambda}_{t,\delta,F+G} = \inf_{|\xi|=1} \sum_{r=1}^{\infty} \int_t^{t+\delta} \left(\sum_{i=1}^{q} \xi_i D_s^{1,r}(F^i+G^i)\right)^2 ds$$

$$\geq \frac{1}{2} \inf_{|\xi|=1} \sum_{r=1}^{\infty} \int_t^{t+\delta} \left(\sum_{i=1}^{q} \xi_i D_s^{1,r} F^i\right)^2 ds$$

$$- \sup_{|\xi|=1} \sum_{r=1}^{\infty} \int_t^{t+\delta} \left(\sum_{i=1}^{q} \xi_i D_s^{1,r} G^i\right)^2 ds$$

$$= \frac{1}{2} \underline{\lambda}_{t,\delta,F} - \overline{\lambda}_{t,\delta,G}. \qquad \square$$

2.2. *Short-time behavior and density evaluations.* We consider some measurable processes $h^{ij}(s), s \in [t, t+\delta]$, $i = 1, \ldots, q$, $j \in N$, such that $h^{ij}(s)$ is $F_t$-measurable and we assume that $\sum_{i=1}^{q} \sum_{j=1}^{\infty} \int_t^{t+\delta} |h^{ij}(s)|^2 ds < \infty$. We define

$$J^i(h) = \sum_{j=1}^{\infty} \int_t^{t+\delta} h^{ij}(s) \, dB^j(s).$$

Since $h(s)$ is $F_t$-measurable, conditionally with respect to $F_t$, $J(h)$ is a Gaussian vector with covariance matrix

$$C^{ij}(J(h)) = \sum_{k=1}^{\infty} \int_t^{t+\delta} h^{ik}(s) h^{jk}(s) \, ds.$$

Given some $F_t$-measurable random variable $V = (V^1, \ldots, V^q)$, we define $G = V + J(h)$. That is,

$$G^i = V^i + \sum_{j=1}^{\infty} \int_t^{t+\delta} h^{ij}(s) \, dB^j(s) = V^i + J^i(h), \qquad i = 1, \ldots, q.$$

Moreover, we consider a deterministic symmetric positive definite matrix $M$ and denote by $\Delta_M$ the smallest eigenvalue of $M$. We assume that $\Delta_M > 0$ (so that $M$ is invertible) and define

$$\|x\|_{M^{-1}} = \sqrt{\langle M^{-1}x, x \rangle},$$

where $\langle \cdot, \cdot \rangle$ is the scalar product on $R^q$.

Given a point $z \in R^q$, a number $a \geq 1$ and a set $A \subseteq \{\omega : \|V(\omega) - z\|_{M^{-1}} \leq 1\}$ we consider the following hypothesis. For every $\omega \in A$,

$(H_1, a, A, z) \qquad aM \geq C(J(h))(\omega) \geq M.$

Note that, in particular, we have

$$\det M \leq \det C(J(h))(\omega) \leq a^q \det M.$$

Finally, we consider $\phi \colon R^q \to R$ defined by $\phi(x) = c \exp(-\frac{1}{1-\|x\|^2})$ for $\|x\| < 1$ and $\phi(x) = 0$ for $\|x\| \geq 1$, with $c$ chosen such that $\int \phi = 1$. We construct the sequence $\phi_\eta \to \delta_0$ defined by $\phi_\eta(y) = \eta^{-q} \phi(\eta^{-1} y)$.

LEMMA 5. *Suppose that $(H_1, a, A, z)$ holds true. Then for every $\eta \in (0, \sqrt{\Delta_M})$,*

$$E_t(\phi_\eta(G-z))(\omega) \geq \frac{1}{e^2 (2\pi a)^{q/2} \sqrt{\det M}} \qquad \text{for } \omega \in A.$$

PROOF. Conditionally with respect to $F_t$, $G - z$ is a Gaussian random variable and so we have

$$E_t(\phi_\eta(G-z)) = \int \phi_\eta(y) \frac{1}{(2\pi)^{q/2} \sqrt{\det C(J(h))}}$$
$$\times \exp\left(-\frac{1}{2}\|y - (V-z)\|^2_{C(J(h))^{-1}}\right) dy.$$

If $\phi_\eta(y) \neq 0$, then $\|y\|_{C(J(h))^{-1}} \leq \|y\|/\sqrt{\Delta_M} \leq \eta/\sqrt{\Delta_M} \leq 1$. Moreover, in view of our hypothesis, $\|V - z\|_{C(J(h))^{-1}} \leq \|V - z\|_{M^{-1}} \leq 1$ so that

$$\exp(-\tfrac{1}{2}\|y - (V-z)\|^2_{C(J(h))^{-1}}) \geq e^{-2}.$$

Since $\int \phi_\eta = 1$ and $\det C(J(h)) \leq a^q \det M$, the proof is completed. □

The following evaluation concerns a perturbation of $G$ by means of a remainder $R$ which is small in an appropriate sense. We consider a $q$-dimensional random variable $R = (R^1, \ldots, R^q)$ such that $R^i \in \bigcap_{p=1}^\infty D^{q+2,p}$, $i = 1, \ldots, q$, and we define

$$F = G + R = V + J(h) + R = V + M^{1/2}(I(h) + R_M)$$

with $R_M = M^{-1/2} R$ and $I(h) = M^{-1/2} J(h) = J(M^{-1/2} h)$. Here, $M^{1/2}$ is also a symmetric invertible positive definite matrix such that $M^{1/2} M^{1/2} = M$. The random variable $I(h)$ will play a role in the following reasoning. Although not a standard normal random variable, it is not far from this; more precisely, under the hypothesis $(H_1, a, A, z)$, one has

(7) $\qquad 1 \leq \inf_{\|\xi\|=1} \langle C(I(h))\xi, \xi\rangle \leq \sup_{\|\xi\|=1} \langle C(I(h))\xi, \xi\rangle \leq a \qquad$ on $A$.

Given $z, A, a$ as in $(H_2, a, A, z)$, we assume that for every $\omega \in A$,

$(H_2, a, A, z) \qquad \|R_M\|_{t,\delta,q+2,p_q} \leq \dfrac{1}{a^{4(q+1)^2} C_q}$

with

(8)
$$p_q = 2^{2(q+2)}p^*(q),$$
$$C_q = c^*(q)\mu(q+1)e^2(2\pi)^{q/2}4^{3(q+3)^3}(q+1)^{q+3},$$

where $c^*(q), p^*(q)$ are those which occur in Proposition 3 and $\mu(q)$ is that which occurs in (2).

REMARK 6. The above constants are neither optimal nor particularly important. What we need is simply to express $C_q$ and $p_q$ as universal constants depending only on the dimension $q$.

REMARK 7. In concrete applications we verify that $\|R_M\|_{t,\delta,q+3,p} \leq C\Delta^\lambda$ for some $\lambda, \Delta > 0$. We then require that $\Delta^\lambda$ be sufficiently small in order to verify the above inequality. In the context of diffusion processes, $\lambda = \frac{1}{2}$ and in the context of the stochastic heat equation, $\lambda = \frac{1}{4}$ (see [13] or [2]).

We also let

$$\Lambda^2 := |DR_M|^2_{t,\delta,1} = \int_t^{t+\delta} \sum_{i=1}^q \sum_{l=1}^\infty |D_s^{1,l} R_M^i|^2 \, ds$$

and note that

$$\overline{\lambda}_{t,\delta,R_M} = \sup_{|\xi|=1} \sum_{i,j=1}^q \xi_i \xi_j \langle DR_M^i, DR_M^j \rangle_{t,\delta,1}$$
$$= \sup_{|\xi|=1} \int_t^{t+\delta} \sum_{l=1}^\infty \langle \xi, D_s^{1,l} R_M \rangle^2 \, ds \leq \Lambda^2.$$

The key evaluation in our approach is given by the following proposition.

PROPOSITION 8. *We consider a point $z \in R^q$, a set $A \subseteq \{\omega : \|V(\omega) - z\|_{M^{-1}} \leq 1\}$ and some $\eta \in (0, \sqrt{\Delta_M})$. Suppose that ($H_1, a, A, z$) and ($H_1, a, A, z$) hold true. Then*

(9)
$$p_\eta(z)(\omega) := E_t(\phi_\eta(F-z))(\omega)$$
$$\geq \frac{1}{4e^2(2\pi a)^{q/2}\sqrt{\det M}} \qquad \text{for } \omega \in A.$$

REMARK 9. Let us give the main ideas of the proof. We write $\Omega = \Gamma \cup \Gamma^c$, where $\Gamma$ is a set on which the Malliavin covariance matrix of $F$ is sufficiently large and $\Gamma^c$ is a set on which we do not control the nondegeneracy of $F$, but which is supposed to be small [in the proof below, we work with $\Theta$ which is a regularization of the indicator function of $\Gamma$, so $\Omega = \Gamma \cup \Gamma^c$ becomes $1 = \Theta + (1 - \Theta)$]. The key question is how we control things on $\Gamma^c$, when we have no information about the nondegeneracy. We write

$$E_t \phi_\eta(F - z) = E_t(\phi_\eta(F - z)\mathbb{1}_\Gamma) + E_t(\phi_\eta(F - z)\mathbb{1}_{\Gamma^c})$$
$$\geq E_t(\phi_\eta(F - z)\mathbb{1}_\Gamma).$$

The above inequality allows us to ignore $\Gamma^c$. Then, following an idea in [13] we use a development in Taylor series of order one to obtain $E_t(\phi_\eta(F - z)\mathbb{1}_\Gamma) = E_t(\phi_\eta(G - z)\mathbb{1}_\Gamma) + r$, where $r$ is a reminder. We are able to evaluate this remainder using Malliavin's integration by parts formula because we are working on the set $\Gamma$ on which the Malliavin covariance matrix is under control. In order to control the first term, we have to "come back to the whole space," that is, to write $E_t(\phi_\eta(G - z)\mathbb{1}_\Gamma) = E_t(\phi_\eta(G - z)) - E_t(\phi_\eta(G - z)\mathbb{1}_{\Gamma^c})$. The previous lemma gives the needed evaluations for $E_t(\phi_\eta(G - z))$, which is the principal term, but we must also evaluate $E_t(\phi_\eta(G - z)\mathbb{1}_{\Gamma^c})$. But (contrary to $F$), $G$ is nondegenerate on the whole space, so we are able to treat this term, even if we are on $\Gamma^c$.

PROOF OF PROPOSITION 8. Since $t$ and $\delta$ are fixed, we will drop them from the notation. So, we write $\phi_F$ instead of $\phi_{t,\delta,F}$, $\|F\|_{k,p}$ instead of $\|F\|_{t,\delta,k,p}$, and so on.

*Step* 1. *Localization.* In the sequel, we will assume (without special mention) that we are on the set $A$ and, in particular, that $\|V - z\|_{M^{-1}} \leq 1$. Since $I(h)$ is Gaussian, its Malliavin covariance matrix coincides with the usual covariance matrix. Moreover, $\underline{\lambda}_{I(h)} \geq 1$ [see (7)] and so, using Proposition 4, we obtain, for every $\rho \in (0, 1)$,

$$(\det \phi_{I(h) + \rho R_M})^{1/q} \geq \tfrac{1}{2}\underline{\lambda}_{I(h)} - \overline{\lambda}_{\rho R_M} \geq \tfrac{1}{2} - \rho\overline{\lambda}_{R_M} \geq \tfrac{1}{2} - \Lambda^2.$$

The aim of this step is to localize on the set on which $\Lambda \leq 1/2$ and consequently $\det \phi_{I(h) + \rho R_M} \geq 1/4^q$. We consider a localization function $\theta \in C_b^\infty(R_+; R_+)$ such that $0 \leq \theta \leq 1$, $\theta(x) = 1$ if $x < 1/4$ and $\theta(x) = 0$ if $x > 1/2$ and we denote $\Theta = \theta(\Lambda)$. We may (see [10], Chapter 1) choose $\theta$ such that for every $0 \leq k \leq q + 1$, one has $\|\theta^{(k)}\|_\infty \leq m(q) := 4^{q+3}(q+1)^{q+3}$.

*Step* 2. *Sobolev norms.* Let us evaluate the Sobolev norm of $\Theta$. First, it is clear that $\|\Theta\|_p = (E_t|\Theta|^p)^{1/p} \leq 1$. Using (2) and (23), we obtain

$$\|\theta(\Lambda)\|_{q+1,2p^*} \leq \mu(q+1)m(q)\|\Lambda\|_{q+1,2^{q+2}p^*}$$
$$\leq 2^{q+1}\mu(q+1)m(q)\|R\|_{q+2,2^{2q+3}p^*} \leq 1,$$

the last inequality being a consequence of $(H_2, a, A, z)$. We conclude that

$$\|\Theta\|_{q+1,2p^*} \leq 2.$$

We now evaluate the Sobolev norms of $1-\Theta$. Since $0 \leq \Theta \leq 1$ and $\Theta = 1$ on the set defined by $\Lambda \leq 1/4$, we have

$$\|1-\Theta\|_p \leq P_t\left(\Lambda \geq \frac{1}{4}\right)^{1/p} = P_t\left(\Lambda^p \geq \frac{1}{4^p}\right)^{1/p}$$

$$\leq 4E_t(\Lambda^p)^{1/p} \leq 16\|\|R_M\|\|_{1,2p}^2.$$

Since $D^k(1-\Theta) = -D^k\Theta$, the same evaluation as for $\Theta$ gives $\|\|1-\Theta\|\|_{q+1,p^*} \leq 2^{q+1}\mu(q+1)m(q)\|\|R\|\|_{q+2,2^{2q+2}p^*}$ so that

$$\|1-\Theta\|_{q,p^*} \leq (16 + 2^{q+1}\mu(q+1)m(q))\|\|R\|\|_{q+2,2^{2q+2}p^*}.$$

Finally, we evaluate the norm of $I(h)$. Note that $D_s^j I^i(h) = (M^{-1/2}h_s)^{ij}$, $i = 1,\ldots,q$, $j \in N$ and $D^\alpha I(h) = 0$ for $|\alpha| > 1$. Since $\int_t^{t+\delta} |\langle D_s I^i(h), D_s I^j(h)\rangle_1|\,ds = (M^{-1/2}C(J(h))(M^{-1/2})^*)^{ij} \leq a$, we use (7) and obtain $\|\|I(h)\|\|_{k,p} \leq a$.

*Step* 3. *Development in Taylor series of order one.* We first localize (multiply by $\Theta$) and then use a development in Taylor series with respect to $R$ in order to obtain

$$p_\eta(z) \geq E_t(\phi_\eta(V-z+J(h)+R)\Theta)$$

$$= E_t(\phi_\eta(G-z)\Theta) + \int_0^1 E_t(\langle\nabla\phi_\eta(V-z+J(h)+\rho R), R\rangle\Theta)\,d\rho$$

$$=: A + B.$$

Let us now evaluate the remainder $B$. We define

$$\Phi_\eta(x) =: \int_{-\infty}^{x_1} dy_1 \cdots \int_{-\infty}^{x_q} dy_q\, \phi_\eta(V-z+M^{1/2}y)$$

so that

$$\phi_\eta(V-z+M^{1/2}x) = \frac{\partial^q \Phi_\eta}{\partial x_1 \cdots \partial x_q}(x).$$

We also have

$$\nabla\phi_\eta(V-z+M^{1/2}x) = M^{-1/2}\nabla\left(\frac{\partial^q \Phi_\eta}{\partial x_1 \cdots \partial x_q}\right)(x) \quad \text{and} \quad 0 \leq \Phi_\eta(x) \leq \frac{1}{\sqrt{\det M}}.$$

The last inequality is obtained using the substitution $x = V - z + M^{1/2}y$ and the fact that $\int \phi_\eta = 1$.

We now write

$$E_t(\langle\nabla\phi_\eta(V-z+J(h)+\rho R), R\rangle\Theta)$$

$$= E_t(\langle\nabla\phi_\eta(V-z+M^{1/2}(I(h)+\rho R_M)), R\rangle\Theta)$$

$$= E_t\left(\left\langle\nabla\left(\frac{\partial^q \Phi_\eta}{\partial x_1 \cdots \partial x_q}\right)(I(h)+\rho R_M), R_M\right\rangle\Theta\right).$$

We use Malliavin's integration by parts formula $q+1$ times to obtain

$$\left|E_t\left(\frac{\partial^{q+1}\Phi_\eta}{\partial x_i\,\partial x_1\cdots\partial x_q}(I(h)+\rho R_M)R_M^i\Theta\right)\right|$$
$$=|E_t(\Phi_\eta(I(h)+\rho R_M)H_{(1,2,\ldots,q,i)}(I(h)+\rho R_M,R_M^i\Theta))|$$
$$\leq \frac{1}{\sqrt{\det M}}E_t(|H_{(1,2,\ldots,q,i)}(I(h)+\rho R_M,R_M^i\Theta)|).$$

We will use Proposition 3 in order to evaluate the above term. Using the notation from Theorem 2 with $F=I(h)+\rho R_M, G=\Theta R_M$ and $k=q+1$, we define the localization set $A:=\{\Lambda\leq 1/4\}$ [note that, since $\theta^{(i)}(\Lambda)=0$, $i=0,\ldots,q+1$, for $\Lambda>1/4$, $R_M^i\Theta$ and its Malliavin derivatives vanish on $A^c$]. Since $\det\phi_{I(h)+\rho R_M}\geq 1/4^q$ on $A$, we obtain

$$v_{p^*}(I(h)+\rho R_M,A)=(E_t((1\vee\det\phi_{I(h)+\rho R_M})^{-p^*}\mathbb{1}_{\{\Lambda\leq 1/4\}}))^{1/p^*}\leq 4^q.$$

Moreover, using (22) and the evaluations from step 2,

$$\|R_M\Theta\|_{q+1,p^*}\leq 2^{q+1}\|R_M\|_{q+1,2p^*}\|\Theta\|_{q+1,2p^*}\leq 2^{q+2}\|R_M\|_{q+1,2p^*}.$$

Finally, using the evaluations for $I(h)$ and hypothesis $(H_2,a,A,z)$, we obtain

$$\|I(h)+\rho R_M\|_{q+2,p^*}\leq \|I(h)\|_{q+2,p^*}+\|R_M\|_{q+2,p^*}\leq 2a.$$

It follows (see Proposition 3) that

$$E_t(|H_{(1,2,\ldots,q,i)}(I(h)+\rho R_M,R_M^i\Theta)|)$$
$$\leq c^*\times 4^{q(q+1)^2}\times 2^{q+2}\|R_M\|_{q+1,2p^*}\times (2a)^{(q+1)(4q-1)}$$
$$\leq c^*\times 4^{3(q+1)^3}\|R_M\|_{q+1,2p^*}\times a^{(q+1)(4q-1)}$$

and, consequently, that

$$|E_t(\langle\nabla\phi_\eta(V-z+J(h)+\rho R),R\rangle\Theta)|$$
$$\leq q\times\frac{a^{(4q-1)(q+1)}}{\sqrt{\det M}}\times c^*\times 4^{3(q+1)^3}\|R_M\|_{q+1,2p^*}$$
$$\leq \frac{1}{4e^2(2\pi a)^{q/2}\sqrt{\det M}},$$

the last inequality being a consequence of $(H_2,a,A,z)$. Finally, $B$ satisfies the same inequality.

We now evaluate $A$. We use the previous lemma to obtain

$$E_t(\phi_\eta(G-z))\geq \frac{1}{e^2(2\pi a)^{q/2}\sqrt{\det M}}.$$

In order to evaluate $E_t(\phi_\eta(G-z)(1-\Theta))$, we integrate by parts $q$ times with respect to $I(h)$ obtain

$$\begin{aligned}E_t(\phi_\eta(G-z)(1-\Theta)) &= E_t(\phi_\eta(V-z+M^{1/2}I(h))(1-\Theta))\\ &= E_t\left(\frac{\partial^q \Phi_\eta}{\partial x_1 \cdots \partial x_q}(I(h))(1-\Theta)\right)\\ &= E_t(\Phi_\eta(I(h))H_{(1,2,\ldots,q)}(I(h),(1-\Theta)))\\ &\leq (\det M)^{-1/2}E_t(|H_{(1,2,\ldots,q)}(I(h),(1-\Theta))|).\end{aligned}$$

Using Proposition 3, the fact that $\underline{\lambda}_{I(h),t} \geq 1$ and the evaluations from step 2, we see that the above term is dominated by

$$\begin{aligned}&(\det M)^{-1/2}c^* \times \|1-\Theta\|_{q,p^*} \times \|I(h)\|_{q+2,p^*}^{q(4(q-1)-1)}\\ &\leq (\det M)^{-1/2}c^* \times (16+2^{q+1}\mu(q+1)m(q))\|R\|_{q+2,2^{2q+2}p^*} \times a^{q(4q-5)}\\ &\leq \frac{1}{2e^2(2\pi a)^{q/2}\sqrt{\det M}},\end{aligned}$$

the last inequality being a consequence of the hypothesis $(H_2, a, A, z)$. It follows that

$$E_t(\phi_\eta(G-z)\Theta) = E_t(\phi_\eta(G-z)) - E_t(\phi_\eta(G-z)(1-\Theta))$$
$$\geq \frac{1}{2e^2(2\pi a)^{q/2}\sqrt{\det M}}.$$

Finally, using the evaluation of $|B|$, we obtain

$$p_\eta(z) \geq \frac{1}{2e^2(2\pi a)^{q/2}\sqrt{\det M}} - |B| \geq \frac{1}{4e^2(2\pi a)^{q/2}\sqrt{\det M}}$$

and the proof is completed. $\square$

2.3. *Evolution sequences.* In this section, the following objects are given:
- A time grid $\Pi_N = (t_0, \ldots, t_N)$ with $0 = t_0 < t_1 < \cdots < t_N = T$. We denote $\delta_k = t_k - t_{k-1}$.
- A sequence of matrixes $M_k$, $k = 0, \ldots, N$, which are deterministic, symmetric, positive definite and invertible. We define by $\Delta_k$ the lower eigenvalue of $M_k$ and define the norms

$$\|x\|_k = \|x\|_{M_k^{-1}} = \sqrt{\langle M_k^{-1}x, x\rangle}.$$

Clearly, $\|x\|_k \leq \Delta_k^{-1/2}\|x\|$. We also consider a sequence of numbers $H_k \geq 1$ such that $H_k^2 M_k \geq M_{k-1}$ in the matrix sense. This is equivalent to

$$\|x\|_k \leq H_k \|x\|_{k-1}$$

and this is the inequality that we need.
- A sequence of numbers $a_k \geq 1$, $k = 0, \ldots, N$.
- A sequence of points $x_k \in R^q$, $k = 1, \ldots, N$, such that

$$\|x_k - x_{k-1}\|_k \leq \tfrac{1}{4}.$$

- A sequence of measurable processes $h_k^{ij}(s)$, $s \in [t_{k-1}, t_k]$, $i = 1, \ldots, q$, $j \in N$, such that $h_k^{ij}(s)$ is $F_{t_{k-1}}$-measurable and $\sum_{i=1}^{q} \sum_{j=1}^{\infty} \int_{t_{k-1}}^{t_k} |h_k^{ij}(s)|^2 \, ds < \infty$ almost surely. We define

$$J_k^i = \sum_{j=1}^{\infty} \int_{t_{k-1}}^{t_k} h_k^{ij}(s) \, dB^j(s).$$

Conditionally with respect to $F_{t_{k-1}}$, $J_k$ is a Gaussian vector with covariance matrix

$$C^{ij}(J_k) = \sum_{l=1}^{\infty} \int_{t_{k-1}}^{t_k} h_k^{il}(s) h_k^{jl}(s) \, ds.$$

- We now introduce the main object of this section, the evolution sequences. We consider a sequence of $R^q$-valued random variables $F_0, \ldots, F_N$ of the form

$$F_k = F_{k-1} + \sum_{j=1}^{\infty} \int_{t_{k-1}}^{t_k} h_k^j(s) \, dB^j(s) + R_k = F_{k-1} + J_k + R_k,$$

where $R_k$ are $q$-dimensional $F_{t_k}$-measurable random variables. In particular, $F_0$ is a constant.

We are interested in the density of the conditional law of $F_k$ with respect to $F_{t_{k-1}}$. Since we do not know that a conditional density exists, we work with the following "regularization of the conditional density":

$$p_{\eta,k}(z) = E_{t_{k-1}}(\phi_\eta(F_k - z)).$$

This quantity makes sense independently of any nondegeneracy assumption.

- Finally, we define the sets

(10) $$A_k = \{\omega : \|F_{i-1}(\omega) - x_i\|_i < \tfrac{1}{2}, i = 1, \ldots, k\} \in F_{t_{k-1}}.$$

DEFINITION 10. We say that $F_0, \ldots, F_N$ is an *elliptic evolution sequence* if $R_k^i \in \bigcap_{p=1}^{\infty} D^{q+2,p}$, $i = 1, \ldots, q$, $k = 1, \ldots, N$, and, on the set $A_k$, one has

(H, i) $$a_k M_k \geq C(J_k) \geq M_k,$$

(H, ii) $$\|M_k^{-1/2} R_k\|_{t_{k-1}, \delta_k, q+2, p_q} \leq \frac{1}{C_q a_k^{4(q+1)^2}},$$

where $C_q$ and $p_q$ are as given in (8).

REMARK 11. The framework of the above definition is inspired by the one introduced by Kohatsu-Higa in [13].

The time grid $\Pi_N$, the path $x = (x_0, \ldots, x_N)$ and $M_k, \Delta_k, a_k, H_k, k = 1, \ldots, N$, are the parameters of evolution sequence and all evaluations are given in terms of these parameters.

As a consequence of the proposition from the previous section, we have the following result:

PROPOSITION 12. *Let $F_0, \ldots, F_N$ be an elliptic evolution sequence and let $k \in \{1, \ldots, N\}$. For every $z \in R^q$ such that $\|x_k - z\|_k \leq \frac{1}{2}$ and every $0 < \eta \leq \sqrt{\Delta_k}$, one has*

$$p_{\eta,k}(z) \geq \frac{1}{4e^2(2\pi a_k)^{q/2}\sqrt{\det M_k}} \quad \text{on the set } A_k.$$

PROOF. Suppose that we are on the set $A_k$. Since $\|x_k - z\|_k \leq \frac{1}{2}$, we have $\|F_{k-1} - z\|_k \leq \|F_{k-1} - x_k\|_k + \|x_k - z\|_k \leq \frac{1}{2} + \frac{1}{2} = 1$ and so $A_k \subseteq \{\|F_{k-1} - z\|_k \leq 1\}$. Since we have an elliptic sequence, the hypothesis $(H_1, a_k, A_k, z)$ and $(H_2, a_k, A_k, z)$ hold true and we may employ Proposition 12. □

2.4. *Tubes evaluations.* The aim of this section is to give lower bounds for $P(A_N)$. We first prove the following lemma:

LEMMA 13. *For every $\eta \in (0, \frac{1}{4H_k}\sqrt{\Delta_{k-1}})$, one has*

$$(11) \quad P(A_k) \geq E\bigg(\mathbb{1}_{A_{k-1}} \int_{\{\|y - x_{k-1}\|_{k-1} \leq (1/4H_k) - \eta/\sqrt{\Delta_{k-1}}\}} p_{\eta,k-1}(y)\, dy\bigg).$$

PROOF. We write

$$P(A_k) = E(\mathbb{1}_{A_{k-1}} E_{t_{k-2}}(\mathbb{1}_{\{\|F_{k-1} - x_k\|_k \leq 1/2\}}))$$

$$= E\bigg(\mathbb{1}_{A_{k-1}} E_{t_{k-2}}\bigg(\int dy\, \phi_\eta(F_{k-1} - y) \mathbb{1}_{\{\|F_{k-1} - x_k\|_k \leq 1/2\}}\bigg)\bigg),$$

the second equality being a consequence of $\int \phi_\eta(F_{k-1} - y)\, dy = \int \phi_\eta(y)\, dy = 1$.

Using the time–space relation $\|x_{k-1} - x_k\|_k \leq 1/4$ and the definition of $H_k$, we obtain

$$\|F_{k-1} - x_k\|_k \leq \tfrac{1}{4} + \|F_{k-1} - x_{k-1}\|_k \leq \tfrac{1}{4} + H_k \|F_{k-1} - x_{k-1}\|_{k-1}$$
$$\leq \tfrac{1}{4} + H_k(\|F_{k-1} - y\|_{k-1} + \|y - x_{k-1}\|_{k-1}).$$

If $\phi_\eta(F_{k-1} - y) \neq 0$, then $\|F_{k-1} - y\| \leq \eta$ and so $\|F_{k-1} - y\|_{k-1} \leq \eta/\sqrt{\Delta_{k-1}}$. Consequently,

$$\|F_{k-1} - x_k\|_k \leq \frac{1}{4} + H_k\left(\frac{\eta}{\sqrt{\Delta_{k-1}}} + \|y - x_{k-1}\|_{k-1}\right).$$

Moreover, if $\|y - x_{k-1}\|_{k-1} \leq (1/4H_k) - \eta/\sqrt{\Delta_{k-1}}$, then $\|F_{k-1} - x_k\|_k \leq 1/2$ and so we may drop this restriction from the integral. We obtain

$$E_{t_{k-2}}\left(\int dy\, \phi_\eta(F_{k-1} - y) \mathbb{1}_{\{\|F_{k-1} - x_k\|_k \leq 1/2\}}\right)$$

$$\geq \int_{\{\|y - x_{k-1}\|_{k-1} \leq (1/4H_k) - \eta/\sqrt{\Delta_{k-1}}\}} E_{t_{k-2}}(\phi_\eta(F_{k-1} - y))\, dy$$

and the proof is completed. □

COROLLARY 14. *Let $F_k, k = 0, \ldots, N$, be an elliptic evolution sequence. For every $k = 1, \ldots, N$,*

(12) $$P(A_k) \geq \frac{1}{8^{q+1} H_k^q e^2 (2qa_{k-1}\pi)^{q/2}} P(A_{k-1}).$$

*In particular,*

(13) $$P(A_N) \geq \left(\frac{1}{8^{q+1} e^2 (2q\pi)^{q/2}}\right)^{N-1} \prod_{k=1}^{N-1} \frac{1}{a_k^{q/2} H_{k+1}^q} \geq e^{-Nq\theta}$$

*with*

(14) $$\theta = \ln(8^2 e^2 (2q\pi)^{1/2}) + \frac{1}{2N} \sum_{k=1}^{N-1} \ln a_k + \frac{1}{N} \sum_{k=2}^{N} \ln H_k.$$

PROOF. We take $\eta = \frac{1}{8H_k}\sqrt{\Delta_{k-1}}$ so that $(1/4H_k) - \eta/\sqrt{\Delta_{k-1}} = 1/8H_k$. Using Lemma 13,

$$P(A_k) \geq E\left(\mathbb{1}_{A_{k-1}} \int_{\{\|y - x_{k-1}\|_{k-1} \leq 1/8H_k\}} p_{\eta,k-1}(y)\, dy\right).$$

Since $H_k \geq 1$, one has $\|y - x_{k-1}\|_{k-1} \leq 1/8H_k \leq 1/2$. We are on the set $A_{k-1}$ so we obtain $\|F_{k-2} - x_{k-1}\|_{k-1} \leq \|F_{k-2} - x_{k-1}\|_{k-1} + \|y - x_{k-1}\|_{k-1} \leq 1$. So, we may use Proposition 12 in order to obtain a lower bound for $p_{\eta,k-1}(y)$ and then

$$P(A_k) \geq \frac{1}{4e^2(2\pi a_{k-1})^{q/2}\sqrt{\det M_{k-1}}} m\left(\|y - x_{k-1}\|_{k-1} \leq \frac{1}{8H_k}\right) P(A_{k-1}),$$

where $m$ is the Lebesgue measure. We use a change of variable and the inequality $m(\|x\| \leq r) \geq (r/\sqrt{q})^q$ in order to obtain

$$\frac{1}{\sqrt{\det M_{k-1}}} m\left(\|y - x_{k-1}\|_{k-1} \leq \frac{1}{8H_k}\right) \geq (8H_k\sqrt{q})^{-q}.$$

It follows that

$$P(A_k) \geq \frac{1}{4e^2(2\pi a_{k-1})^{q/2}} \times \frac{1}{q^{q/2}8^q H_k^q} \times P(A_{k-1})$$

and (12) is proved. In order to prove (13), we employ recurrence to obtain

$$P(A_N) \geq \left(\frac{1}{8^{q+1}e^2(2q\pi)^{q/2}}\right)^{N-1} \prod_{k=2}^{N} \frac{1}{a_{k-1}^{q/2} H_k^q} P\left(\|F_0 - x_1\|_1 \leq \frac{1}{2}\right).$$

Since $\|F_0 - x_1\|_1 = \|x_0 - x_1\|_1 \leq \frac{1}{2}$, (13) is proved. □

2.5. *The main result.* Our final result is as follows. We look for lower bounds for the density of $F_N$. We say that the law of $F_N$ has a local density $p_{F_N}$ in a neighborhood of $x_N$ with respect to the Lebesgue measure on $R^q$ if there exists some $\delta > 0$ such that for every smooth function $\psi$ with the support included in the ball $B_\delta(x_N)$, one has

$$E\psi(F_N) = \int \psi(x) p_{F_N}(x) \, dx.$$

THEOREM 15. *Let $F_k, k = 0, \ldots, N$ be an elliptic evolution sequence. Suppose that the law of $F_N$ has a continuous local density $p_{F_N}$ in a neighborhood of $x_N$ with respect to the Lebesgue measure on $R^q$. Then*

$$p_{F_N}(x_N) \geq \frac{1}{4e^2(2\pi a_N)^{q/2}\sqrt{\det M_N}} e^{-Nq\theta}$$

*with $\theta$ given as in (14).*

PROOF. We use Proposition 12 and the fact that $A_N$ is $F_{t_{N-1}}$-measurable to obtain

$$\int_{R^q} p_{F_N}(x)\phi_\eta(x - x_N) \, dx = E(\phi_\eta(F_N - x_N)) = E(E_{t_{N-1}}\phi_\eta(F_N - x_N))$$

$$\geq E(E_{t_{N-1}}(\phi_\eta(F_N - x_N))\mathbb{1}_{A_N})$$

$$\geq \frac{1}{4e^2(2\pi a_N)^{q/2}\sqrt{\det M_N}} P(A_N)$$

$$\geq \frac{1}{4e^2(2\pi a_N)^{q/2}\sqrt{\det M_N}} e^{-Nq\theta}.$$

We now use the continuity of $p_{F_N}$ and take the limit with $\eta \to 0$ in order to obtain the result. □

**3. Elliptic Itô processes.** We consider a $q$-dimensional Itô process of the form

$$X_t^i = x_0^i + \sum_{j=1}^\infty \int_0^t U_s^{ij}\, dB_s^j + \int_0^t V_s^i\, ds, \qquad i=1,\ldots,q,$$

and assume that for every $T > 0$,

(i) $\quad E\left(\int_0^T (\|U_s\|^2 + |V_s|)\, ds\right) < \infty,$

(ii) $\quad U_s, V_s \in \bigcap_{p \in N} D^{q+2,p} \qquad \forall 0 \leq s \leq T$

with

$$\|U_s\|^2 = \sum_{i=1}^q \sum_{j=1}^\infty |U_s^{ij}|^2 \quad \text{and} \quad |V_s|^2 = \sum_{i=1}^q |V_s^i|^2.$$

We fix $T > 0$ and $y \in R^q$ and study the density of the law of $X_T$ in $y$. In order to do this, we have to give a nondegeneracy assumption on $X_T$ and this assumption is related to a deterministic path from $x_0$ to $y$, that is, a continuously differentiable function $x:[0,T] \to R^q$ such that $x(0) = x_0$ and $x(T) = y$. We also consider continuous, strictly positive functions $r, K:[0,T] \to R_+$ and a number $a \geq 1$. The significance of these functions is as follows. We work on a tube around the deterministic path $x(t)$; $r(t)$ represents the radius of this tube and $K$ controls the small increments of our process. The number $a$ comes on in the ellipticity assumption. Finally, we consider a family $Q_t$, $t \in [0,T]$, of symmetric, positive definite and invertible matrixes. We denote

$$C^{ij}(U_t) = \sum_{l=1}^\infty U^{il}(t) U^{jl}(t),$$

$$\Gamma_\delta^i(t) = \sum_{j=1}^\infty \int_t^{t+\delta} (U_s - U_t)^{ij}\, dB_s^j + \int_t^{t+\delta} V_s^i\, ds.$$

Our ellipticity hypothesis as follows.

DEFINITION 16. Let $\nu > 0$. We say that the path $x$ is $(r,K,a,Q,\nu)$-elliptic for $X$ if for every $0 < t < T$ and $0 < \delta < T - t$,

(H, i) $\qquad\qquad a \times Q_t \geq C(U_t) \geq Q_t,$

(H$^\nu$, ii) $\qquad \|Q_t^{-1/2} \Gamma_\delta(t)\|_{t,\delta,q+2,p_q} \leq K(t)\delta^{1/2} + \nu$

on the set $\{\omega : \|X(t,\omega) - x(t)\|_{Q^{-1}(t)} \leq r(t)\}$. Recall that $p_q$ is given in (8).

We need some more notation. Given $m \geq 1, h > 0$ we denote by $L(m, h)$ the class of the strictly positive functions $f : [0, T] \to R_+$ such that

(15) $$f(s) \leq m f(t) \qquad \text{for } |s - t| \leq h.$$

If the above inequality holds true for every $t$ and $s$ in $[0, T]$, then we take $h = \infty$.

Moreover, we say that $Q \in L(m, h)$ if

(16) $$\|x\|_{Q_s^{-1}} \leq m \|x\|_{Q_t^{-1}} \qquad \text{for } t \leq s \leq t + h.$$

Note that this is still equivalent to $Q_t \leq m^2 Q_s$.

THEOREM 17. *We suppose that the law of $X_T$ has a continuous local density in $y$ and that there exists a path $x(t)$, $t \in [0, T]$ such that $x(0) = x, x(T) = y$ and which is $(K, a, r, Q, \nu)$-elliptic for $X$. We also consider two functions $\pi_t, \gamma_t$ such that*

$$\pi(t) \leq \min\{r^2(t), (C_q K_t a^{4(q+1)^2})^{-1/\nu}\},$$

$$\|\partial_t x_t\|_{Q_t^{-1}} = \sqrt{\langle Q_t^{-1} \partial_t x_t, \partial_t x_t \rangle} \leq \gamma(t) \qquad \forall t \in [0, T],$$

*where $C_q$ is given in (8). Finally we assume that $Q \in L(m_Q, h_Q)$, $\pi \in L(m_\pi, h_\pi)$, $\gamma \in L(m_\gamma, h_\gamma)$ for some constants $m_Q, m_\pi, m_\gamma \geq 1$ and $h_Q, h_\pi, h_\gamma > 0$. We denote $h := h_Q \wedge h_\pi \wedge h_\gamma$. Then*

$$p_T(x_0, y) \geq \frac{1}{4e^2 (2\pi T m_Q a)^{q/2} \sqrt{\det Q_T}}$$

$$\times \exp\left(-q\left(\alpha + \frac{a}{2}\right) \times \int_0^T \left(\frac{1}{h} + \frac{m_\pi}{\pi(t)} + 16 m_Q^2 \gamma_t^2\right) dt\right),$$

*where $\alpha = \ln(8e(2\pi q)^{1/4}) + \ln m_Q + 4 \ln m_\gamma + \ln m_\pi$.*

REMARK 18. We may take $\gamma_t = \|\partial_t x_t\|_{Q_t^{-1}}$, but in concrete examples, it may be difficult to work with this function (to compute $m_\gamma$, e.g.)—this is why we allow $\gamma_t$ to be larger. The same holds for $\pi$.

PROOF OF THEOREM 17. *Step* 1. We define a time grid $t_k, k \in N$, in the following way. We take $t_0 = 0$ and, if $t_k$ is given, we define

$$\tau_k = \inf\left\{u > 0 : \int_{t_k}^{t_k + u} \gamma_s^2 \, ds \geq \frac{1}{16 m_Q^2}\right\}, \qquad t_{k+1} = t_k + h \wedge \pi(t_k) \wedge \tau_k.$$

We put $N = \min\{k : t_k \geq T\}$ and claim that

(17) $$N \leq \int_0^T \left(\frac{1}{h} + \frac{m_\pi}{\pi(t)} + 16 m_Q^2 \gamma_t^2\right) dt.$$

In order to prove this, we denote $I = \{k \leq N : t_{k+1} - t_k = \tau_k\}$, $I' = \{k \leq N : t_{k+1} - t_k = \pi(t_k)\}$ and $I'' = \{k \leq N : t_{k+1} - t_k = h\}$ and write

$$\int_0^T \left(\frac{1}{h} + \frac{m_\pi}{\pi(t)} + 16 m_Q^2 \gamma_t^2\right) dt \geq \sum_{k \in I} \int_{t_k}^{t_k + \tau_k} 16 m_Q^2 \gamma_t^2 \, dt$$

$$+ \sum_{k \in I'} \int_{t_k}^{t_k + \pi(t_k)} \frac{m_\pi}{\pi(t)} \, dt + \sum_{k \in I''} \int_{t_{k-1}}^{t_k + h} \frac{1}{h} \, dt.$$

We claim that all terms in the above sums are greater than one; hence, (17) holds true. For $k \in I$, this follows from the definition of $\tau_k$ and for $k \in I''$, it is trivial. Suppose, now, that $k \in I'$ and note that in this case, $\pi(t_k) \leq h \leq h_\pi$. Then $\pi(t) \leq m_\pi \pi(t_k)$ for $t_k \leq t \leq t_k + \pi(t_k)$ and so

$$\int_{t_k}^{t_k + \pi(t_k)} \frac{m_\pi}{\pi(t)} \, dt \geq 1.$$

The proof of (17) is thus completed.

*Step* 2. We define an evolution sequence as follows. We define $\delta_k := t_k - t_{k-1}$, $F_k = X(t_k)$ and write

$$F_k = F_{k-1} + \sum_{j=1}^\infty \int_{t_{k-1}}^{t_k} U_s^j dB_s^j + \int_{t_{k-1}}^{t_k} V_s \, ds = F_{k-1} + J_k + R_k$$

with $J_k = \sum_{j=1}^\infty \int_{t_{k-1}}^{t_k} U^j(t_{k-1}) \, dB_s^j$ and $R_k = \Gamma_{\delta_k}(t_{k-1})$.

Coming back to the notation concerning the evolution sequences, we have $h_k(s) = U(t_{k-1})$ for $s \in [t_{k-1}, t_k)$ and so $C(J_k) = \delta_k C(U(t_{k-1}))$, $M_k = \delta_k Q(t_{k-1})$ and $a_k = a$.

*Step* 3. We denote $H_k^2 =: m_Q(m_\gamma^4 \vee m_\pi)$ and check that $H_k^2 M_k \geq M_{k-1}$, that is, $H_k^2 \delta_k Q(t_{k-1}) \geq \delta_{k-1} Q(t_{k-2})$. Since $\delta_k \leq h \leq h_Q$, we use (16) to obtain $m_Q Q(t_{k-1}) \geq Q(t_{k-2})$. So, it suffices to show that $(H_k^2/m_Q)\delta_k \geq \delta_{k-1}$, which reads

$$(*) \qquad H(h \wedge \pi(t_{k-1}) \wedge \tau_{k-1}) \geq h \wedge \pi(t_{k-2}) \wedge \tau_{k-2}$$

$$\text{with } H := \frac{H_k^2}{m_Q} = m_\gamma^4 \vee m_\pi.$$

Since $H \geq m_\pi$ and $t_{k-1} - t_{k-2} \leq h \leq h_\pi$, we have $H\pi(t_{k-1}) \geq \pi(t_{k-2})$ and so

$$(**) \qquad H(h \wedge \pi(t_{k-1})) \geq h \wedge \pi(t_{k-2}).$$

If $H\tau_{k-1} > h \geq h \wedge \pi_h(t_{k-2}) \wedge \tau_{k-2}$, then $(**)$ guarantees that $(*)$ holds true. We now consider the case where $H\tau_{k-1} \leq h$ and, in particular, $\tau_{k-1} \leq h \leq h_\gamma$ (because $H \geq 1$). For every $t \in [t_{k-2}, t_{k-2} + H\tau_{k-1})$, one has $|t - t_{k-1}| \leq h \leq h_\gamma$ so we may use (15) (twice) to obtain

$$\int_{t_{k-2}}^{t_{k-2}+H\tau_{k-1}} \gamma_t^2 \, dt \geq \frac{1}{m_\gamma^2} \gamma_{t_{k-1}}^2 H\tau_{k-1}$$

$$\geq \frac{H}{m_\gamma^4} \int_{t_{k-1}}^{t_{k-1}+\tau_{k-1}} \gamma_t^2 \, dt = \frac{H}{m_\gamma^4} \frac{1}{16 m_Q^2} \geq \frac{1}{16 m_Q^2}.$$

This proves that $H\tau_{k-1} \geq \tau_{k-2}$ and this, together with (**), gives (*).

*Step* 4. Our aim is to check that $F_k, k = 0, \ldots, N$, is an elliptic evolution sequence (see Definition 10). We take $x_k =: x(t_k)$. We will first check the space–time relation $\|x(t_k) - x(t_{k-1})\|_k \leq \frac{1}{4}$. We write $x(t_k) - x(t_{k-1}) = \int_{t_{k-1}}^{t_k} \partial_t x_t \, dt$ and using (16), we obtain

$$\|x(t_k) - x(t_{k-1})\|_k = \frac{1}{\sqrt{\delta_k}} \left\| \int_{t_{k-1}}^{t_k} \partial_t x_t \, dt \right\|_{Q^{-1}(t_{k-1})} \leq \frac{1}{\sqrt{\delta_k}} \int_{t_{k-1}}^{t_k} \|\partial_t x_t\|_{Q^{-1}(t_{k-1})} \, dt$$

$$\leq \frac{m_Q}{\sqrt{\delta_k}} \int_{t_{k-1}}^{t_k} \|\partial_t x_t\|_{Q^{-1}(t)} \, dt \leq \frac{m_Q}{\sqrt{\delta_k}} \int_{t_{k-1}}^{t_k} \gamma(t) \, dt$$

$$\leq m_Q \left( \int_{t_{k-1}}^{t_k} \gamma^2(t) \, dt \right)^{1/2} \leq \frac{1}{4},$$

the last inequality being a consequence of the definition of $\tau_k$. So, the space–time relation is verified.

Moreover, if $\omega \in A_k$ [see (10)], we have $\|x(t_k) - X(t_{k-1})\|_k \leq \frac{1}{2}$ and so

$$\|x(t_{k-1}) - X(t_{k-1})\|_{Q^{-1}(t_{k-1})}$$
$$\leq \sqrt{\delta_k}(\|x(t_k) - x(t_{k-1})\|_k + \|x(t_k) - X(t_{k-1})\|_k)$$
$$\leq \sqrt{\delta_k}(\tfrac{1}{4} + \tfrac{1}{2}) \leq \sqrt{\delta_k} \leq r(t_{k-1}),$$

the last inequality being true because $\delta_k \leq \pi(t_{k-1}) \leq r(t_{k-1})^2$. We have proved that $A_k \subseteq \{\|x(t_{k-1}) - X(t_{k-1})\|_{Q^{-1}(t_{k-1})} \leq r(t_{k-1})\}$, so we may use the hypotheses (H, i), (H$^\nu$, ii).

Using (H, i), we obtain $a_k M_k = a\delta_k Q(t_{k-1}) \geq C(J_k) = \delta_k C(U(t_{k-1})) \geq \delta_k Q(t_{k-1}) = M_k$. In particular, $\det C(J_k) \leq a\delta_k \det Q(t_{k-1})$.

Since $M_k^{-1/2} R_k = \delta_k^{-1/2} Q^{-1/2}(t_{k-1}) \Gamma_{\delta_k}(t_{k-1})$, the hypothesis (H$^\nu$, ii) gives

$$\|M_k^{-1/2} R_k\|_{t_{k-1}, \delta_k, q+2, p_q} \leq K(t_{k-1}) \delta_k^\nu \leq \frac{1}{C_q a^{4(q+1)^2}},$$

the last inequality being a consequence of $\delta_k \leq \pi(t_{k-1})$. So, we have proved that we have an elliptic evolution sequence.

*Step* 5. We are now able to use the density evaluations from the previous section. We note that $\frac{1}{N}\sum_{k=1}^{N} \ln H_k \leq \ln m_Q + 4\ln m_\gamma + \ln m_\pi$ and $\frac{1}{2N}\sum_{k=1}^{N} \ln a = \frac{a}{2}$, there by obtaining

$$\theta := \frac{1}{2}\ln(8^2 e^2 (2\pi q)^{1/2}) + \frac{1}{2qN}\sum_{k=1}^{N} \ln a + \frac{1}{N}\sum_{k=2}^{N} \ln H_k \leq \alpha + \frac{a}{2}.$$

Further, note also that $a^{q/2}(\det M_N)^{1/2} = a^{q/2}\delta_N^{q/2}\sqrt{\det Q(t_{N-1})} \leq a^{q/2}m_Q^{q/2} \times T^{q/2}\sqrt{\det Q(T)}$. Finally, we use Theorem 17 and our evaluation (17) for $N$ to obtain

$$p_T(x_0, y) \geq \frac{1}{4e^2(2\pi a)^{q/2}\sqrt{\det M_N}}e^{-Nq\theta}$$

$$\geq \frac{1}{4e^2(2\pi Tm_Q a)^{q/2}\sqrt{\det Q(T)}}$$

$$\times \exp\left(-q\left(\alpha + \frac{a}{2}\right) \times \int_0^T \left(\frac{1}{h} + \frac{e_\pi(h)}{\pi(t)} + 16m_Q^2\gamma_t^2\right)dt\right)$$

and the proof is thus completed. □

**4. Diffusion processes.** In this section, we will study the diffusion process $X$ which is the solution of the SDE

$$dX_t^i = \sum_{j=1}^{d} \sigma_j^i(X_t)\,dB_t^j + b^i(X_t)\,dt, \qquad i=1,\ldots,q,$$

$$X_0 = x_0.$$

We fix $\varepsilon_i \in \{0,1\}$, $i = 0,\ldots,q$, and denote

$$N^2(x) =: \varepsilon_0 + \sum_{i=1}^{q} \varepsilon_i |x^i|^2.$$

We assume that the coefficients $\sigma$ and $b$ are of class $C^{q+2}$ and verify that

(A, i) $\qquad\qquad\max_i(\|\sigma^i(x)\| + |b^i(x)|) \leq C_0 N(x),$

(A, ii) $\qquad\max_i(\|\sigma^i(x) - \sigma^i(y)\| + |b^i(x) - b^i(y)|) \leq C_0\|x-y\|,$

(A, iii) $\qquad\max_{|\alpha|\leq q+2}\max_{i,j}(|D^\alpha\sigma_j^i(x)| + |D^\alpha b^i(x)|) \leq C_0.$

The reason to use $N(x)$ (instead of the usual Euclidean norm) in order to control the growth of the coefficients is that for different choices of $\varepsilon_i$, $i = 0,\ldots,q$, we obtain different type of hypothesis—bonded coefficients, linear growth, log-normal types diffusions, and so on—and the behavior of the lower bound of the density is different in these cases.

As an immediate consequence of (A), one has

(A, iv) $$\lambda^*(x) = \sup_{\|\xi\|=1} \langle \sigma\sigma^*(x)\xi, \xi\rangle \leq qC_0^2 N^2(x),$$

(A, v) $$\sup_{\|\xi\|=1} |\langle \sigma\sigma^*(x)\xi, \xi\rangle - \langle \sigma\sigma^*(y)\xi, \xi\rangle| \leq qC_0^2(2N(x) + \|x-y\|)\|x-y\|,$$

(A, vi) $$|\det \sigma\sigma^*(x) - \det \sigma\sigma^*(y)| \leq q! C_0^{2q}(2N(x) + \|x-y\|)^{2q-1}\|x-y\|.$$

It is clear that $X$ is an Itô process and that

$$\Gamma^i_\delta(t) = \sum_{j=1}^d \int_t^{t+\delta} (\sigma^{ij}(X_s) - \sigma^{ij}(X_t))\,dB^j_s + \int_t^{t+\delta} b^i(X_s)\,ds.$$

We will employ the following standard lemma:

LEMMA 19. *Suppose that* (A) *holds true. Then for every* $t \geq 0$, $1 \geq \delta > 0$, $0 \leq m \leq q+2$, $p \in N$, *one has*

(18) $$\|\Gamma_\delta(t)\|_{t,\delta,m,p} \leq C(m,p)N(X_t)\delta,$$

*where* $C(m,p)$ *is a constant which depends on* $C_0$ *and on* $m,p$.

PROOF. The proof is straightforward, but rather long and tedious, so we just outline the main arguments (see [11] for a complete approach to such evaluations). In order to simplify the notation, we take $b=0$. The first step is to check that for $t \leq s \leq \delta \leq 1$, $(E_t(|N(X_s)|^p))^{1/p} \leq CN(X_t)$. Here and in the sequel, $C$ is a constant which may change from one line to another. We use the *SDE* of $X$, Hölder's inequality, Burckholder's inequality and hypothesis (A, i) in order to obtain $E_t(|X^i_s|^p) \leq C|X^i_t|^p + C'\int_t^s E_t|N(X_r)|^p\,dr$. We then take $\sum_{i=1}^q \varepsilon_i E_t(|X^i_s|^p)$ and employ Gronwell's lemma. This proves the above inequality. The same argument gives

$$E_t(|X^i_s - X^i_t|^p) \leq C\delta^{p/2-1}\int_t^s E_t|N(X_r)|^p\,dr$$
$$\leq CN^p(X_t)\delta^{p/2}.$$

It follows that

$$E_t(|\Gamma^i_\delta(s)|^p) \leq C\sum_{j=1}^q E_t\left(\int_t^{t+\delta}|\sigma^i_j(X_r) - \sigma^i_j(X_t)|^2\,dr\right)^{p/2}$$
$$\leq CN^p(X_t)\delta^p.$$

Let us now deal with the first-order Malliavin derivatives. For $t \leq u \leq s \leq t' \leq t + \delta$ and $i = 1, \ldots, q$, $l = 1, \ldots, d$, one has

$$D_u^l X_s^i = \sigma_l^i(X_u) + \sum_{j=1}^d \int_u^s \nabla \sigma_j^i(X_r) D_u^l X_r \, dB_r^j.$$

We look to $s \to D_u X_s^i = (D_u^1 X_s^i, \ldots, D_u^d X_s^i)$ as an $R^d$-valued process and use Burckholder's inequality (for $R^d$-valued martingales) in order to obtain

$$E_t \|D_u X_s^i\|^p \leq C E_t \|\sigma^i(X_u)\|^p + C \sum_{j=1}^d E_t \left\| \int_u^s \nabla \sigma_j^i(X_r) D_u^l X_r^i \, dB_r^j \right\|^p$$

$$\leq C E_t |N(X_u)|^p + C \sum_{j=1}^d E_t \left( \int_u^s \|\nabla \sigma_j^i(X_r) D_u X_r\|^2 \, dr \right)^{p/2}$$

$$\leq C N(X_t)^p + C E_t \left( \int_u^s \|D_u X_r\|^2 \, dr \right)^{p/2}.$$

It follows that

$$E_t \left( \int_u^{t'} \|D_u X_s^i\|^2 \, ds \right)^{p/2}$$

$$\leq C \delta^{p/2-1} \int_u^{t'} E_t \|D_u X_s^i\|^p \, ds$$

$$\leq C \delta^{p/2} N(X_t)^p + C \delta^{p/2-1} \int_u^{t'} E_t \left( \int_u^s \|D_u X_r\|^2 \, dr \right)^{p/2} ds.$$

Using Gromwell's lemma, we obtain

$$E_t \left( \int_u^{t'} \|D_u X_s^i\|^2 \, ds \right)^{p/2} \leq C \delta^{p/2} N(X_t)^p.$$

Finally, for $u \in [t, t + \delta]$, one has $D_u^l \Gamma_\delta^i(t) = \sigma_l^i(X_u) - \sigma_l^i(X_t) + \sum_{j=1}^d \int_u^{t+\delta} \nabla \sigma_j^i(X_r) D_u^l X_r^i \, dB_r^j$ and, so, using (A, ii), Hölder's inequality and Burckholder's inequality we obtain

$$E_t \left( \int_t^{t+\delta} \|D_u \Gamma_\delta^i(t)\|^2 du \right)^{p/2}$$

$$\leq C \delta^{p/2-1} \int_t^{t+\delta} E_t \|D_u \Gamma_\delta^i(t)\|^p \, du$$

$$\leq C \delta^{p/2-1} \int_t^{t+\delta} E_t \|\sigma^i(X_u) - \sigma^i(X_t)\|^p \, du$$

$$+ C \delta^{p/2-1} \int_t^{t+\delta} \sum_{j=1}^d E_t \left\| \int_u^{t+\delta} \nabla \sigma_j^i(X_r) D_u^l X_r^i \, dB_r^j \right\|^p du$$

$$\leq C\delta^{p/2-1} \int_t^{t+\delta} E_t \|X_u - X_t\|^p \, du$$

$$+ C\delta^{p/2-1} \int_t^{t+\delta} E_t \left( \int_u^{t+\delta} \|D_u^l X_r^i\|^2 \, dr \right)^{p/2} du$$

$$\leq C\delta^p N(X_t)^p.$$

So, we have proved that $\|\Gamma_\delta(t)\|_{t,\delta,1,p} \leq CN(X_t)\delta$. The proof is analogous for higher order derivatives, so we omit it. $\square$

We denote by $\lambda_*$ the smallest eigenvalue of $\sigma\sigma^*$ and let

$$\rho(x) := \frac{\sqrt{\lambda_*(x)}}{N(x)}.$$

Roughly speaking, $\rho^2$ is of the same order as the quotient of the smallest and the largest eigenvalues of $\sigma\sigma^*$.

LEMMA 20. (i) *Suppose that* (A) *holds true and let* $x = (x_t)_{t \leq T}$ *be a differentiable path such that* $\lambda_*(x_t) > 0$, *for all* $0 \leq t \leq T$. *Then the path $x$ is* $(r, K, a, Q, \nu)$-*elliptic (in the sense of Definition* 16*) with* $\nu = \frac{1}{2}, a = 3/2$ *and*

(19)
$$Q_t = \frac{1}{2}\sigma\sigma^*(x_t), \qquad r_t = \frac{\rho^2(x_t)}{6q^{3/2}C_0^3},$$

$$K_t = C(q+2, p_q)\left(\frac{1}{\rho(x_t)} + \frac{1}{\sqrt{\lambda_*(x_t)}}\right),$$

*where* $C(q+2, p_q)$ *is the constant from* (18) *and* $p_q$ *is given in* (8).

(ii) *Assume that there exists a measurable function* $M_t$, $t \in [0, T]$, *and a number* $h_G \in (0, 1)$ *such that for every* $t \in [0, T]$,

(G)
$$\|\partial x_t\| \leq M_t N(x_t),$$

$$h_G \int_t^{t+h_G} M_s^2 \, ds \leq \frac{1}{4q}.$$

*Then for every* $s, t \in [0, T]$ *such that* $|s - t| \leq h_G$, *we have*

$$N(x_s) \leq 4N(x_t).$$

PROOF. Suppose that $\|X_t - x_t\|_{Q_t^{-1}} \leq r_t$. In view of (A, iv), $\lambda^*(x) \leq qC_0^2 \times N^2(x)$ and so $(\sigma\sigma^*)^{-1}(x) \geq (1/qC_0^2 N^2(x)) \times I$, where $I$ is the identity matrix. It follows that $\|X_t - x_t\| \leq \sqrt{q}C_0 N(x_t)\|X_t - x_t\|_{Q_t^{-1}} \leq \sqrt{q}C_0 N(x_t)r_t$. Let $\xi \in R^q$ with $\|\xi\| = 1$. Using (A, v) and $\sqrt{q}C_0 r_t \leq 1$, we obtain

$$|\langle \sigma\sigma^*(X_t)\xi,\xi\rangle - \langle \sigma\sigma^*(x_t)\xi,\xi\rangle| \leq qC_0^2(2N(x_t) + \sqrt{q}C_0N(x_t)r_t)\sqrt{q}C_0N(x_t)r_t$$
$$\leq 3q^{3/2}C_0^3 N^2(x_t)r_t \leq \frac{\lambda_*(x_t)}{2},$$

the last inequality being a consequence of the choice of $r_t$. This gives

$$\langle \sigma\sigma^*(X_t)\xi,\xi\rangle = \langle \sigma\sigma^*(x_t)\xi,\xi\rangle + (\langle \sigma\sigma^*(X_t)\xi,\xi\rangle - \langle \sigma\sigma^*(x_t)\xi,\xi\rangle)$$
$$\geq \langle \sigma\sigma^*(x_t)\xi,\xi\rangle - \frac{\lambda_*(x_t)}{2} \geq \frac{1}{2}\langle \sigma\sigma^*(x_t)\xi,\xi\rangle = \langle Q_t\xi,\xi\rangle.$$

Moreover,

$$\langle \sigma\sigma^*(X_t)\xi,\xi\rangle \leq \langle \sigma\sigma^*(x_t)\xi,\xi\rangle + \frac{\lambda_*(x_t)}{2} \leq \frac{3}{2}\langle Q_t\xi,\xi\rangle.$$

So (H, i) holds true with $a_k = 3/2$.

Note that $N(X_t) \leq N(x_t) + r_t \leq N(x_t) + 1$. Using the previous lemma with $m = q+2$ and $p = p_q$, we obtain

$$\|Q_t^{-1}\Gamma_\delta(t)\|_{t,\delta,m,p} \leq \frac{1}{\sqrt{\lambda_*(x_t)}}\|\Gamma_\delta(t)\|_{t,\delta,m,p} \leq \frac{C(m,p)}{\sqrt{\lambda_*(x_t)}}N(X_t)\delta$$
$$\leq \frac{C(m,p)}{\sqrt{\lambda_*(x_t)}}(N(x_t)+1)\delta = K_t\delta.$$

So, we have an elliptic path with parameters given in (19).

Let us now prove (ii). Suppose that $t < s$ and write $x_s = x_t + \int_t^s \partial x_r\, dr$ so that

$$N^2(x_s) = \varepsilon_0 + \sum_{i=1}^q \varepsilon_i |x_s^i|^2 \leq \varepsilon_0 + 2\sum_{i=1}^q \varepsilon_i |x_t^i|^2 + 2\sum_{i=1}^q \varepsilon_i \left|\int_t^s \partial x_r^i\, dr\right|^2$$
$$\leq 2N^2(x_t) + 2\sum_{i=1}^q \varepsilon_i (s-t)\int_t^s |M_r N(x_r)|^2\, dr.$$

By the choice of $h_G$,

$$\sup_{t\leq s\leq t+h_G} N^2(x_s) \leq 2N^2(x_t) + 2q \sup_{t\leq s\leq t+h_G} N^2(x_s) h_G \int_t^{t+h_G} M_r^2\, dr$$
$$\leq 2N^2(x_t) + \tfrac{1}{2}\sup_{t\leq s\leq t+h_G} N^2(x_s)$$

and the proof is completed. $\square$

We are now able to state our result.

THEOREM 21. *Suppose that* (A) *holds true and that* $x = (x_t)_{t \leq T}$ *is a differentiable path such that* $x_0 = x_0$, $x_T = y$ *and* $\rho(x_t) \geq \frac{1}{\mu}$, $\lambda_*^{-2}(x_t) \leq \chi$ *for all* $0 \leq t \leq T$, *for some* $\mu \geq 1, \chi > 0$. *We assume that there exists a number* $h_G \in [0,1)$ *and a measurable function* $M$ *such that* (G) *holds true and* $M \in L(\eta_M, h_M)$ *for some* $\eta_M \geq 1, h_M > 0$. *Then*

$$p_T(x_0, y) \geq \frac{1}{4e^2(6\mu\sqrt{q}\pi T)^{q/2}\sqrt{\det \sigma\sigma^*(y)}}$$

(20)
$$\times \exp\Bigl(-K_q T(1 + \ln C_0 + \ln \mu + \ln \eta_M)$$

$$\times \Bigl(\frac{C_0^2 \mu^4}{T}\int_0^T M_r^2\, dr$$

$$+ \mu^4 \vee (\mu+\chi)^2 K_{\text{diff}} + \frac{1}{h_G \wedge h_M}\Bigr)\Bigr)$$

*with* $K_{\text{diff}} = C_0^6 C^2(q+2, p_q)$ [*recall that* $C_0$ *is given in hypothesis* (A) *and* $C(q+2, p_q)$ *given in* (18)] *and* $K_q$ *is a constant depending only on* $q$.

REMARK 22. Usually, the constants which appear in the lower bound are independent of $x_0, y$ and $T$, but the dependence on the coefficients $\sigma, b$ is not explicit. So, the lower bound is not significant for $y$ in a compact set, but only for $y \to \infty$. Here, the constants are explicit (although not optimal), so the result is relevant for every $x_0, y$—this is the motivation of the (rather tedious) effort to keep the constants under control.

PROOF OF THEOREM 21. Under our assumptions, $\sigma\sigma^*(x_0) > 0$ and so the law of $X_T$ has a continuous density with respect to the Lebesgue measure. We will use Theorem 2 in order to obtain the lower bound. By the previous lemma, $x$ is $(r, K, a, Q, \nu)$-elliptic and we know the corresponding parameters [see (19)]. Since $\rho^{-1}(x_t) \leq \mu$, we take

$$K_t = C(q+2, p_q)(\mu + \chi), \qquad r_t = \frac{1}{6\mu^2 q^{3/2} C_0^3}.$$

We also have

$$r_t^2 \wedge \frac{1}{C_q^2 K_t^2 a_t^{8(q+1)^2}} = \frac{1}{36\mu^4 q^3 C_0^6}$$

$$\wedge \Bigl(\frac{1}{C_q^2 (3/2)^{8(q+1)^2}} \times \frac{1}{(\mu+\chi)^2 C^2(q+2, p_q)}\Bigr)$$

$$\geq \frac{1}{K_q \times K_{\text{diff}}} \times \frac{1}{\mu^4 \vee (\mu+\chi)^2} =: \pi_t.$$

Since the function $\pi$ is constant, $m_\pi = 1$ and $h_\pi = \infty$. Moreover, using (G), we have

$$\|\partial x_t\|^2_{Q_t^{-1}} \leq \frac{2}{\lambda_*(x_t)}\|\partial x_t\|^2 \leq \frac{2N^2(x_t)}{\lambda_*(x_t)}M_t^2 \leq 2\mu^2 M_t^2,$$

so we take $\gamma_t = \sqrt{2}\mu M_t$ and have $m_\gamma = \eta_M$ and $h_\gamma = h_M$.

We now take $h_Q = h_G$ and compute $m_Q$. Using point (ii) for the previous lemma, $N^2(x_t)/N^2(x_s) \leq 16$ for $|s - t| \leq h_G$. Moreover, by (A, iv),

$$\frac{2}{qC_0^2 N^2(x_t)} \leq Q_t^{-1} \leq \frac{2}{\lambda_*(x_t)}$$

so that

$$Q_s^{-1} \leq \frac{2}{\lambda_*(x_s)} = \frac{qC_0^2}{\rho^2(x_s)} \times \frac{N^2(x_t)}{N^2(x_s)} \times \frac{2}{qC_0^2 N^2(x_t)} \leq 16qC_0^2 \theta^2 Q_t^{-1}.$$

So, we take $m_Q = 4\sqrt{q}C_0\mu$. Finally, $h = h_Q \wedge h_\pi \wedge h_\gamma = h_G \wedge h_M$. We compute

$$\alpha = \ln(8e(2q\pi)^{1/4}) + \ln m_Q + 4\ln m_\gamma + \ln m_\pi$$
$$= \ln(32qe(2q\pi)^{1/4}) + \ln C_0 + \ln\mu + 4\ln\eta_M.$$

We now use the evaluation from Theorem 2 to obtain

$$p_T(x_0, y) \geq \frac{1}{4e^2(6\mu\sqrt{q}\pi T)^{q/2}\sqrt{\det\sigma\sigma^*(y)}}$$
$$\times \exp\Big(-K_q(1 + \ln C_0 + \ln\mu + \ln\eta_M)$$
$$\times \int_0^T \Big(\frac{1}{h_M \wedge h_G} + \mu^4 \vee (\mu + \chi)^2 K_{\text{diff}} + C_0^2\mu^4 M_t^2\Big)dt\Big)$$

and the proof is completed. $\square$

We will now write the path $x$ in a special form, given by (a variant of) the skeleton of the diffusion process $X$. More precisely, we consider some $\phi(t) = (\phi^1(t), \ldots, \phi^d(t))$, $t \in [0, T]$, such that $\phi \in L^2([0,T])^d$ and we denote $\|\phi\|_T^2 = \int_0^T \|\phi_t\|^2 dt$. We associate with $\phi$ the path $x = x^\phi$ which solves the deterministic differential equation

$$(\text{E}_\phi) \qquad dx(t) = \sum_{j=1}^d \sigma_j(x(t))\phi_t^j\, dt, \qquad x(0) = x_0.$$

REMARK 23. Note that for every differentiable path $x_t$ such that $\sigma\sigma^*(x_t) \geq \lambda_*(x_t) > 0$, there exists $\phi$ such that $(\text{E}_\phi)$ holds true. Indeed, if one takes $\phi_t = \sigma^*(x_t)(\sigma\sigma^*)^{-1}(x_t)\partial_t x_t$, then $\sigma(x_t)\phi_t = \partial_t x_t$.

We consider a set of parameters $\theta = (\mu, \chi, \nu, \eta, h)$, $\mu, \nu, \eta \geq 1$, $h, \chi > 0$, and we define $\Phi_\theta(x_0, y)$ to be the set of the controls $\phi \in (L^2([0,T]))^d$ such that

$$x_0^\phi = x_0, \qquad x_T^\phi = y,$$

$$\rho(x_t^\phi) \geq \frac{1}{\mu}, \qquad \sqrt{\lambda_*(x_t^\phi)} \geq \frac{1}{\chi} \qquad \forall t \in [0,T],$$

$$\|\phi_t\| \leq \eta \|\phi_s\| \qquad \forall |s-t| \leq h, \qquad \|\phi\|_t \leq \nu \qquad \forall t \leq T.$$

We then define

$$d_\theta(x_0, y) = \inf\{\|\phi\|_T : \phi \in \Phi_\theta(x_0, y)\}$$
$$= \infty \qquad \text{if } \Phi_\theta(x_0, y) = \varnothing.$$

THEOREM 24. *We assume that (A) holds true. We fix $x_0, y \in R^q$ and $\theta$ and suppose that $d_\theta(x_0, y) < \infty$. Then the law of $X_T$ has a continuous density $p_T(x_0, y)$ which verifies*

$$\begin{aligned}
p_T(x_0, y) &\geq \frac{1}{4e^2(6\mu\sqrt{q}\pi T)^{q/2}\sqrt{\det \sigma\sigma^*(y)}} \\
&\quad \times \exp\Big(-K_q(1 + \ln(C_0\mu\eta)) \\
&\quad \times \Big(C_0^4\mu^4 d_\theta^2(x_0, y) + T\Big(\mu^4 \vee (\mu+\chi)^2 K_{\text{diff}} \\
&\quad\quad + \frac{1}{h} + 2C_0\nu\sqrt{q}\Big)\Big)\Big),
\end{aligned}$$
(21)

*where $K_{\text{diff}} = C_0^6 C^2(q+2, p_q)$ [recall that $C_0$ is given in the hypothesis (A) and $C(q+2, p_q)$ is given in (18)] and $K_q$ is a constant depending only on $q$.*

PROOF. We fix $\phi \in \Phi_\theta$ and take $x$ to be solution of $\partial_t x_t = \sigma(x_t)\phi_t$. Using the orthogonal decomposition $\phi_t = \sigma^*(x_t)v + w$ with $v \in R^d$ and $w$ such that $\sigma(x_t)w = 0$, one obtains $\partial_t x_t = \sigma(x_t)\sigma^*(x_t)v$. Consequently,

$$\langle (\sigma\sigma^*(x_t))^{-1}\partial_t x_t, \partial_t x_t \rangle = \langle v, \sigma\sigma^*(x_t)v \rangle = \|\sigma^*(x_t)v\|^2 \leq \|\phi_t\|^2.$$

Using (A, iv), we obtain

$$\|\partial x_t\|^2 \leq C_0^2 N^2(x_t) < \langle (\sigma\sigma^*(x_t))^{-1}\partial_t x_t, \partial_t x_t \rangle \leq C_0^2 N^2(x_t)\|\phi_t\|^2,$$

so we take $M_t = C_0\|\phi_t\|$ in (G). We take $h_G = 1/2\sqrt{q}C_0\nu$ and these by obtain

$$h_G \int_t^{t+h_G} M_s^2\, ds \leq h_G^2 C_0^2 \nu^2 \leq \frac{1}{4q}.$$

So (G) holds true.

Since $\phi \in \Phi_\theta$, one has $\eta_M = \eta, h_M = h$. Moreover, $\int_0^T M_r^2 \, dr = C_0^2 \|\phi\|_T^2$. Substituting this into (20), we obtain

$$p_T(x_0, y) \geq \frac{2^q}{4e^2(6\mu\sqrt{q}\pi T)^{q/2}\sqrt{\det \sigma\sigma^*(y)}}$$

$$\times \exp\Big(-K_{(q)}q(1 + \ln(C_0\mu\eta))\Big)$$

$$\times \Big(C_0^4\mu^4\|\phi\|_T^2$$

$$+ T\Big(K_{\text{diff}}\mu^4 \vee (\mu + \chi)^2 + 2\sqrt{q}C_0\nu + \frac{1}{h}\Big)\Big).$$

We now take the infimum over $\phi \in \Phi_\theta$ and the proof is completed. □

## APPENDIX

We will use the following Hölder inequalities for the conditional Malliavin norms.

LEMMA 25. *Let* $H, Q \in \bigcap_{p \geq 1} D^{k+1,p}$. *Suppose that* $Q(\omega) = 0$, $D^1 Q(\omega) = 0, \ldots, D^k Q(\omega) = 0$ *for every* $\omega \in A^c$, *where* $A$ *is some measurable set. Then for every* $p \geq 1$,

(22) $$\|HQ\|_{t,\delta,k,p} \leq k! 2^k \|H\|_{t,\delta,k,2p,A} \|Q\|_{t,\delta,k,2p,A},$$

(23) $$\|\langle DH, DQ \rangle_{t,\delta,1}\|_{t,\delta,k,p} \leq k! 2^k \|H\|_{t,\delta,k+1,2p,A} \|Q\|_{t,\delta,k+1,2p,A},$$

*where*

$$\|H\|_{t,\delta,k,p,A}^p := E_t(\mathbb{1}_A \|H\|_{t,\delta,k}^p),$$

$$\|H\|_{t,\delta,k,p,A}^p := \|H\|_{t,\delta,k,p,A}^p - E_t(\mathbb{1}_A |H|^p).$$

PROOF. Let us introduce some notation. Let $I \subseteq \{1, \ldots, k\}$, $I = \{i_1, \ldots, i_r\}$ with $1 \leq i_1 < \cdots < i_r \leq k$. We denote $|I| = r$. Given a multi-index $a = \{a_1, \ldots, a_k\} \in \{1, \ldots, d\}^k$ and $s = (s_1, \ldots, s_k) \in R_+^k$, we denote $a(I) = \{a_{i_1}, \ldots, a_{i_r}\}$ and $s(I) = (s_{i_1}, \ldots, s_{i_r})$. We also put $ds = ds_1 \ldots ds_k$ and $ds(I) = ds_{i_1} \ldots ds_{i_r}$. We denote $D_s^{k,a} = D_{s_k}^{a_k} \ldots D_{s_1}^{a_1}$ and we write

$$D_s^{k,a}(HQ) = \sum_{i=0}^k \sum_{|I|=i} D_{s(I)}^{i,a(I)} H \times D_{s(I^c)}^{k-i,a(I^c)} Q.$$

Since the sum has $2^k$ terms, we have

$$|D_s^{k,a}(HQ)|^2 \leq 2^k \sum_{i=0}^k \sum_{|I|=i} |D_{s(I)}^{i,a(I)} H|^2 \times |D_{s(I^c)}^{k-i,a(I^c)} Q|^2.$$

Let us denote
$$\alpha_a(HQ) := \int_{[t,t+\delta]^k} |D_s^{k,a}(HQ)|^2 \, ds.$$

Since
$$\int_{[t,t+\delta]^k} |D_{s(I)}^{i,a(I)} H \times D_{s(I^c)}^{k-i,a(I^c)} Q|^2 \, ds$$
$$= \int_{[t,t+\delta]^i} |D_{s(I)}^{i,a(I)} H|^2 \, ds(I) \times \int_{[t,t+\delta]^{k-i}} |D_{s(I^c)}^{k-i,a(I^c)} Q|^2 \, ds(I^c)$$
$$= \alpha_{a(I)}(H) \alpha_{a(I^c)}(Q),$$

we obtain

(24) $$\alpha_a(HQ) \leq 2^k \sum_{i=0}^{k} \sum_{|I|=i} \alpha_{a(I)}(H) \alpha_{a(I^c)}(Q).$$

We claim that this implies
$$\|HQ\|_{t,\delta,k}^2 = \sum_{|a|\leq k} \alpha_a(HQ) \leq k! 2^k \sum_{|b|\leq k} \alpha_b(H) \times \sum_{|c|\leq k} \alpha_c(Q)$$
(25)
$$= k! 2^k \|H\|_{t,\delta,k}^2 \|Q\|_{t,\delta,k}^2.$$

In order to check this inequality, we consider two multi-indices, $b$ and $c$, such that $|b| \leq k$ and $|c| \leq k$. The term $\alpha_b(H)\alpha_c(Q)$ will appear in the right-hand side of (24) for any multi-index $a$ such that $b = a(I)$ and $c = a(I^c)$ for some $I$. It follows that the components of $a$ are fixed (they are the reunion of the components of $b$ and $c$). It follows that there are at most $k!$ such terms. So, (25) follows.

Since $Q$ and its derivatives are null on $A^c$, we have $\|H\|_{t,\delta,k}^2 \|Q\|_{t,\delta,k}^2 = \mathbb{1}_A \|H\|_{t,\delta,k}^2 \|Q\|_{t,\delta,k}^2$ and we use Hölder's inequality in order to obtain (22). The proof of (23) is similar. $\square$

PROOF OF LEMMA 1. Let $\Gamma^{ij}$ be the cofactor corresponding to $i,j$ of the matrix $\phi_{t,\delta,F}$ and let $d := \det \phi_{t,\delta,F}$. We have $\widehat{\phi}_{t,\delta,F}^{ij}(\omega) = \frac{1}{d}\Gamma^{ij}(\omega)$ for $\omega \in A$. Using (22),

$$\|(\widehat{\phi}_{t,\delta,F} G)^i\|_{t,\delta,k,p} \leq \sum_{i=1}^{q} \|\widehat{\phi}_{t,\delta,F}^{ij} G^j\|_{t,\delta,k,p} = \sum_{j=1}^{q} \|d^{-1}\Gamma^{ij} G^j\|_{t,\delta,k,p}$$
$$\leq 2^{2(k+1)} \sum_{j=1}^{q} \|d^{-1}\|_{t,\delta,k,4p,A} \|\Gamma^{ij}\|_{t,\delta,k,4p,A} \|G^j\|_{t,\delta,k,4p,A}.$$

Moreover, using (22) and (23), we have $\|\Gamma^{ij}\|_{t,\delta,k,4p} \leq 2^{2(q-1)(k+1)} \times \|F\|_{t,\delta,k+1,2^{2(q+1)}p}^{2(q-1)}$ and $\|d\|_{t,\delta,k,8p} \leq 2^{2q(k+1)}\|F\|_{t,\delta,k+1,2^{2(q+2)}p}^{2q}$. Then using the chain rule for the function $\frac{1}{x}$, we obtain

$$\|d^{-1}\|_{t,\delta,k,4p,A} \leq c(k,q)v_{8(k+1)p}^{k+1}(F,A)\|d\|_{t,\delta,k,8p}$$
$$\leq c(k,q)v_{8(k+1)p}^{k+1}(F,A)\|F\|_{t,\delta,k+1,2^{2(q+2)}p}^{2q},$$

where $c(k,q)$ is generic notation for a constant which depends on $k$ and $q$. It follows that

$$\|(\widehat{\phi}_{t,\delta,F}G)^i\|_{t,\delta,k,p} \leq c(k,q)\|G\|_{t,\delta,k,4p}v_{8(k+1)p}^{k+1}(F,A)\|F\|_{t,\delta,k+1,2^{q+2}p}^{4q-2}$$

and the proof is completed. □

PROOF OF PROPOSITION 3. We denote by $c$ a constant which depends on $k, l$ and $q$ and which may change from one line to the next. Using (22), (23), (4) and (3), we obtain

$$\|H_i(F,G)\|_{t,\delta,k,p} \leq \sum_{j=1}^{q}(\|G\widehat{\phi}_{t,\delta,F}^{ji}L_{t,\delta}(F_j)\|_{t,\delta,k,p}$$
$$+ \|\langle DF_j, D(\widehat{\phi}_{t,\delta,F}^{ji}G)\rangle_{t,\delta,1}\|_{t,\delta,k,p})$$
$$\leq c\sum_{j=1}^{q}(\|G\widehat{\phi}_{t,\delta,F}^{ji}\|_{t,\delta,k,2p}\|L_{t,\delta}(F_j)\|_{t,\delta,k,2p}$$
$$+ \|F_j\|_{t,\delta,k+1,2p}\|G\widehat{\phi}_{t,\delta,F}^{ji}\|_{t,\delta,k+1,2p})$$
$$\leq c\|G\|_{t,\delta,k+1,4p}v_{16(k+1)p}^{k+1}(F,A)\|F\|_{t,\delta,k+2,2^{q+3}p}^{4q-1}.$$

Consider now a multi-index $\alpha$ with $|\alpha|=l$. We iterate the above relation to obtain

$$\|H_\alpha(F,G)\|_{t,\delta,k,p} \leq c\|G\|_{t,\delta,k+l,2^{3l}p}v_{2^{4l}(k+l+1)^lp}^{l(k+l)}(F,A)\|F\|_{t,\delta,k+l+1,2^{l(q+3)}p}^{l(4q-1)}$$

and the proof is completed. □

**Acknowledgments.** I am grateful to Sylvie Meleard and Eulalia Nualart who read an initial version of this paper and pointed out to me an error which has been corrected here. I would also like to thank to Arturo Kohatsu-Higa for several very useful discussions.

## REFERENCES

[1] AIDA, S., KUSUOCKA, S. and STROOCK, D. (1993). On the support of Wiener functionals. In *Asymptotic Problems in Probability Theory*: *Wiener Functionals and Asymptotics* (K. D. Elworthy and N. Ikeda, eds.). *Pitman Res. Notes Math. Ser.* **284** 3–34. Longman Sci. Tech., Harlow. MR1354161

[2] BALLY, V. (2003). Lower bounds for the density of the law of locally elliptic Itô processes. Preprint 4887, INRIA.

[3] BALLY, V. and PARDOUX, E. (1998). Malliavin calculus for white noise driven parabolic SPDEs. *Potential Anal.* **9** 27–64. MR1644120

[4] BASS, R. F. (1997). *Diffusions and Elliptic Operators. Probability and Its Applications.* Springer, New York. MR1483890

[5] BEN-AROUS, G. and LEANDRE, R. (1991). Decroissance exponentille du noyau de la chaleur sur la diagonale. II. *Probab. Theory Related Fields* **90** 377–402. MR1133372

[6] DALANG, R. and NUALART, E. (2004). Potential theory for hyperbolic SPDEs. *Ann. Probab.* **32** 2099–2148. MR2073187

[7] FEFFERMAN, C. L. and SANCHEZ-CALLE, A. (1986). Fundamental solutions of second order subelliptic operators. *Ann. of Math.* (*2*) **124** 247–272. MR0855295

[8] FRIEDMAN, A. (1975). *Stochastic Differential Equations and Applications.* Academic Press, New York. MR0494490

[9] GUÉRIN, H., MÉLÉARD, S. and NULALART, E. (2006). Estimate for the density of a nonlinear Landau proccess. *J. Funct. Anal.* **238** 649–677.

[10] HÖRMANDER, L. (1990). *The Analysis of Linear Partial Differential Operators.* I, 2nd ed. Springer, Berlin. MR1065993

[11] KUSUOKA, S. and STROOCK, D. (1985). Applications of the Malliavin calculus, part II. *J. Fac. Sci. Univ. Tokyo Sect. IA Math.* **32** 1–76. MR0783181

[12] KUSUOKA, S. and STROOCK, D. (1987). Applications of the Malliavin calculus, part III. *J. Fac. Sci. Univ. Tokyo Sect. IA Math.* **34** 391–442. MR0914028

[13] KOHATSU-HIGA, A. (2003). Lower bounds for densities of uniform elliptic random variables on Wiener space. *Probab. Theory Related Fields* **126** 421–457. MR1992500

[14] IKEDA, N. and WATANABE, S. (1989). *Stochastic Differential Equations and Diffusion Processes*, 2nd ed. North-Holland, Amsterdam. MR1011252

[15] MILLET, A. and SANZ, M. (1997). Points of positive density for the solutions to a hyperbolic SPDE. *Potential Anal.* **7** 623–659. MR1473646

[16] NUALART, D. (1995). *The Malliavin Calculus and Related Topics.* Springer, New York. MR1344217

[17] SANCHEZ-CALLE, A. (1986). Fundamental solutions and geometry of the sum of square of vector fields. *Invent. Math.* **78** 143–160. MR0762360

LABORATOIRE D'ANALYSE
 ET DE MATHÉMATIQUES APPLIQUÉES
UMR 8050
UNIVERSITÉ DE MARNE-LA-VALLÉE
CITÉ DESCARTES
5 BLD DESCARTES
CHAMPS-SUR-MARNE
F-77454 MARNE-LA VALLÉE
CEDEX 2
FRANCE
E-MAIL: bally@univ-mlv.fr